\documentclass{amsart}

\usepackage{amsmath}
\usepackage{amsfonts}
\usepackage{amssymb,enumerate}
\usepackage{amsthm}
\usepackage[all]{xy}
\usepackage{hyperref}

\newtheorem{lem}{Lemma}[section]

\newtheorem{prop}[lem]{Proposition}
\newtheorem{thm}[lem]{Theorem}

\newtheorem{Defn}[lem]{Definition}
\newtheorem{Ex}[lem]{Example}
\newtheorem{Question}[lem]{Question}
\newtheorem{Property}[lem]{Property}
\newtheorem{Properties}[lem]{Properties}
\newtheorem{Discussion}[lem]{Remark}
\newtheorem{Construction}[lem]{Construction}
\newtheorem{Subprops}{}[lem]
\newtheorem{Para}[lem]{}
\newtheorem{Notation}[lem]{Notation}
\newtheorem{Fact}[lem]{Fact}
\newtheorem{Notationdefinition}[lem]{Definition/Notation}

\newenvironment{defn}{\begin{Defn}\rm}{\end{Defn}}
\newenvironment{ex}{\begin{Ex}\rm}{\end{Ex}}

\newenvironment{disc}{\begin{Discussion}\rm}{\end{Discussion}}

\newenvironment{notation}{\begin{Notation}\rm}{\end{Notation}}
\newenvironment{fact}{\begin{Fact}\rm}{\end{Fact}}
\newenvironment{notationdefinition}{\begin{Notationdefinition}\rm}{\end{Notationdefinition}}

\newtheorem{intthm}{Theorem}

\newcommand{\cat}[1]{\mathcal{#1}}
\newcommand{\catx}{\cat{X}}
\newcommand{\caty}{\cat{Y}}
\newcommand{\catm}{\cat{M}}
\newcommand{\catv}{\cat{V}}
\newcommand{\catw}{\cat{W}}
\newcommand{\catg}{\cat{G}}
\newcommand{\catp}{\cat{P}}
\newcommand{\catf}{\cat{F}}
\newcommand{\cati}{\cat{I}}
\newcommand{\cata}{\cat{A}}

\newcommand{\catb}{\cat{B}}

\newcommand{\catgic}{\cat{GI}_C}
\newcommand{\catgib}{\cat{GI}_B}
\newcommand{\catgid}{\cat{GI}_D}
\newcommand{\catgicd}{\cat{GI}_{C^{\dagger}}}
\newcommand{\caticd}{\cat{I}_{C^{\dagger}}}
\newcommand{\catgc}{\cat{G}_C}

\newcommand{\catgpc}{\cat{GP}_C}
\newcommand{\catgpb}{\cat{GP}_B}
\newcommand{\catgpd}{\cat{GP}_D}

\newcommand{\catac}{\cat{A}_C}
\newcommand{\catbc}{\cat{B}_C}
\newcommand{\catbb}{\cat{B}_B}
\newcommand{\catacd}{\cat{A}_{\da{C}}}
\newcommand{\catbcd}{\cat{B}_{\da{C}}}
\newcommand{\catgfc}{\cat{GF}_C}
\newcommand{\catic}{\cat{I}_C}
\newcommand{\catib}{\cat{I}_B}
\newcommand{\catpc}{\cat{P}_C}
\newcommand{\catfc}{\cat{F}_C}
\newcommand{\catpb}{\cat{P}_B}
\newcommand{\catid}{\cat{I}_D}
\newcommand{\catpd}{\cat{P}_D}
\newcommand{\PP}{\catp}
\newcommand{\I}{\cati}
\newcommand{\opg}{\cat{G}}
\newcommand{\finrescat}[1]{\operatorname{res}\comp{\cat{#1}}}
\newcommand{\proprescat}[1]{\operatorname{res}\wti{\cat{#1}}}
\newcommand{\finrescatx}{\finrescat{X}}
\newcommand{\finrescaty}{\finrescat{Y}}

\newcommand{\finrescatp}{\finrescat{P}}

\newcommand{\finrescatgpc}{\operatorname{res}\comp{\catgpc}}
\newcommand{\finrescatggpc}{\operatorname{res}\comp{\catg(\catpc)}}

\newcommand{\fincorescatggicd}{\operatorname{cores}\comp{\catg(\caticd)}}
\newcommand{\proprescatpc}{\operatorname{res}\wti{\catpc}}
\newcommand{\finrescatpc}{\operatorname{res}\comp{\catpc}}

\newcommand{\fincorescatic}{\operatorname{cores}\comp{\catic}}
\newcommand{\fincorescaticd}{\operatorname{cores}\comp{\caticd}}
\newcommand{\proprescatgpc}{\operatorname{res}\wti{\catgpc}}

\newcommand{\proprescatx}{\proprescat{X}}

\newcommand{\proprescatw}{\proprescat{W}}
\newcommand{\fincorescat}[1]{\operatorname{cores}\comp{\cat{#1}}}
\newcommand{\propcorescat}[1]{\operatorname{cores}\wti{\cat{#1}}}
\newcommand{\fincorescatx}{\fincorescat{X}}
\newcommand{\fincorescati}{\fincorescat{I}}

\newcommand{\fincorescatgicd}{\operatorname{cores}\comp{\catgicd}}

\newcommand{\propcorescatic}{\operatorname{cores}\wti{\catic}}
\newcommand{\propcorescatgic}{\operatorname{cores}\wti{\catgic}}

\newcommand{\fincorescaty}{\fincorescat{Y}}

\newcommand{\propcorescaty}{\propcorescat{Y}}
\newcommand{\propcorescatv}{\propcorescat{V}}

\newcommand{\G}{\mathcal{G}}
\newcommand{\catggpc}{\opg(\catpc)}

\newcommand{\catggic}{\opg(\catic)}

\newcommand{\pd}{\operatorname{pd}}

\newcommand{\gcdim}[1]{\mathrm{G}_{#1}\text{-}\dim}	
	
\newcommand{\id}{\operatorname{id}}

\newcommand{\Catpd}[1]{\cat{#1}\text{-}\pd}
\newcommand{\xpd}{\Catpd{X}}

\newcommand{\wpd}{\Catpd{W}}

\newcommand{\gpd}{\Catpd{GP}}
\newcommand{\gfd}{\Catpd{GF}}
\newcommand{\gid}{\Catid{GI}}
\newcommand{\Catid}[1]{\cat{#1}\text{-}\id}
\newcommand{\yid}{\Catid{Y}}

\newcommand{\pdpd}{\catpd\text{-}\pd}
\newcommand{\idid}{\catid\text{-}\id}
\newcommand{\pcpd}{\catpc\text{-}\pd}
\newcommand{\pbpd}{\catpb\text{-}\pd}
\newcommand{\icdagdim}{\caticd\text{-}\id}
\newcommand{\icdid}{\caticd\text{-}\id}
\newcommand{\icdim}{\catic\text{-}\id}
\newcommand{\icid}{\catic\text{-}\id}
\newcommand{\ibid}{\catib\text{-}\id}
\newcommand{\pcdim}{\catpc\text{-}\pd}
\newcommand{\gpcpd}{\catgpc\text{-}\pd}
\newcommand{\gfcpd}{\catgfc\text{-}\pd}
\newcommand{\gpbpd}{\catgpb\text{-}\pd}
\newcommand{\gicid}{\catgic\text{-}\id}
\newcommand{\gibid}{\catgib\text{-}\id}
\newcommand{\gicdagdim}{\catgicd\text{-}\id}
\newcommand{\gicdid}{\catgicd\text{-}\id}
\newcommand{\ggpcpd}{\catg(\catpc)\text{-}\pd}

\newcommand{\ggicid}{\catg(\catic)\text{-}\id}

\newcommand{\depth}{\operatorname{depth}}

\newcommand{\len}{\operatorname{length}}

\newcommand{\ext}{\operatorname{Ext}}	
\newcommand{\eext}{\operatorname{e}}	
\newcommand{\weext}{\operatorname{e}_{\catw R}}	
\newcommand{\veext}{\operatorname{e}_{R\catv}}	
\newcommand{\xeext}{\operatorname{e}_{\catx R}}	
\newcommand{\yeext}{\operatorname{e}_{R\caty}}

\newcommand{\HH}{\operatorname{H}}
\newcommand{\Hom}{\operatorname{Hom}}	
\newcommand{\coker}{\operatorname{Coker}}
\newcommand{\spec}{\operatorname{Spec}}
\newcommand{\mspec}{\mathrm{m}\text{-}\spec}

\newcommand{\tor}{\operatorname{Tor}}
\newcommand{\ann}{\operatorname{Ann}}
\newcommand{\im}{\operatorname{Im}}
\newcommand{\shift}{\mathsf{\Sigma}}

\newcommand{\da}[1]{#1^{\dagger}}

\newcommand{\supp}{\operatorname{Supp}}

\newcommand{\cone}{\operatorname{Cone}}

\newcommand{\Ker}{\operatorname{Ker}}

\newcommand{\aext}{\ext_{R}}
\newcommand{\pcext}{\ext_{\catpc}}
\newcommand{\pdext}{\ext_{\catpd}}
\newcommand{\gpdext}{\ext_{\catgpd}}
\newcommand{\pbext}{\ext_{\catpb}}
\newcommand{\gpbext}{\ext_{\catgpb}}
\newcommand{\icext}{\ext_{\catic}}
\newcommand{\ibext}{\ext_{\catib}}
\newcommand{\idext}{\ext_{\catid}}
\newcommand{\gidext}{\ext_{\catgid}}

\newcommand{\icdext}{\ext_{\caticd}}
\newcommand{\gpcext}{\ext_{\catgpc}}
\newcommand{\gicext}{\ext_{\catgic}}
\newcommand{\gibext}{\ext_{\catgib}}

\newcommand{\gicdext}{\ext_{\catgicd}}
\newcommand{\ahom}{\Hom_{R}}
\newcommand{\xaext}{\ext_{\catx R}}
\newcommand{\axext}{\ext_{R\catx}}
\newcommand{\ayext}{\ext_{R\caty}}
\newcommand{\avext}{\ext_{R\catv}}

\newcommand{\waext}{\ext_{\catw R}}
\newcommand{\mat}[1]{#1^{\vee}}

\newcommand{\ideal}[1]{\mathfrak{#1}}
\newcommand{\m}{\ideal{m}}

\newcommand{\p}{\ideal{p}}

\newcommand{\comp}[1]{\widehat{#1}}
\newcommand{\ol}{\overline}
\newcommand{\wti}{\widetilde}

\newcommand{\from}{\leftarrow}
\newcommand{\xra}{\xrightarrow}
\newcommand{\xla}{\xleftarrow}

\newcommand{\xwacomp}{\vartheta_{\catx \catw R}}
\newcommand{\ayvcomp}{\vartheta_{R \caty \catv}}

\newcommand{\xrcomp}{\varkappa_{\catx R}}
\newcommand{\rycomp}{\varkappa_{R\caty}}

\newcommand{\gpcpccomp}{\vartheta_{\catgpc\catpc}}
\newcommand{\gpcpbcomp}{\vartheta_{\catgpc\catpb}}
\newcommand{\gpcgpbcomp}{\vartheta_{\catgpc\catgpb}}
\newcommand{\gicibcomp}{\vartheta_{\catgic\catib}}
\newcommand{\gicgibcomp}{\vartheta_{\catgic\catgib}}
\newcommand{\giciccomp}{\vartheta_{\catgic\catic}}
\newcommand{\pccomp}{\varkappa_{\catpc}}
\newcommand{\gpccomp}{\varkappa_{\catgpc}}
\newcommand{\iccomp}{\varkappa_{\catic}}
\newcommand{\icdcomp}{\varkappa_{\caticd}}
\newcommand{\giccomp}{\varkappa_{\catgic}}

\newcommand{\A}{\mathcal{A}}

\renewcommand{\geq}{\geqslant}
\renewcommand{\leq}{\leqslant}
\renewcommand{\ker}{\Ker}
\renewcommand{\hom}{\Hom}

\begin{document}

\bibliographystyle{amsplain}

\author{Sean Sather-Wagstaff}

\address{Sean Sather-Wagstaff, Department of Mathematical Sciences, Kent State University,
  Mathematics and Computer Science Building, Summit Street, Kent OH
  44242, USA}
\email{sather@math.kent.edu}
\urladdr{http://www.math.kent.edu/~sather}

\author{Tirdad Sharif}
\address{Tirdad Sharif, School of Mathematics, Institute for Studies in
Theoretical Physics and Mathematics, P. O. Box 19395-5746, Tehran, Iran}
\email{sharif@ipm.ir}
\urladdr{http://www.ipm.ac.ir/IPM/people/personalinfo.jsp?PeopleCode=IP0400060}
\thanks{TS is supported by a grant from IPM,(No. 83130311).}

\author{Diana White} 
\address{Diana White, Department of Mathematics, University of Nebraska,
   203 Avery Hall, Lincoln, NE, 68588-0130 USA} \email{dwhite@math.unl.edu}
\urladdr{http://www.math.unl.edu/~s-dwhite14/}

\title[Comparison of relative cohomology theories]{Comparison of relative cohomology theories
with respect to semidualizing modules
}

\dedicatory{Dedicated to the memory of Anders Juel Frankild}

\date{\today}


\begin{abstract}
We compare and contrast various relative cohomology theories that arise
from resolutions involving semidualizing modules.  We prove a general balance result
for relative cohomology over a Cohen-Macaulay ring with a dualizing module,
and we demonstrate the failure of the naive version of balance
one might expect for these functors.  We prove that the natural comparison
morphisms between relative cohomology modules are isomorphisms in
several cases, and we provide a Yoneda-type description
of the first relative Ext functor.
Finally, we show by example that each distinct relative cohomology
construction does in fact result in a different functor.
\end{abstract}

\maketitle

\section*{Introduction} 

The study of relative homological algebra was initiated by
Butler and Horrocks~\cite{butler:cer}  
and 
Eilenberg and Moore~\cite{eilenberg:frha}
and has been revitalized recently by a number of authors,
notably, Enochs and Jenda~\cite{enochs:rha}
and Avramov and Martsinkovsky~\cite{avramov:aratc}.
The basic idea behind this construction is to consider 
resolutions of a module $M$ over a ring $R$, where the modules in the 
resolutions are taken from a fixed class $\catx$. One restricts focus to those
resolutions $X$, called proper $\catx$-resolutions, 
with good enough lifting properties to make them unique
up to homotopy equivalence, 
and this yields well-defined functors
$$\xaext^n(M,-)=\HH_{-n}(\ahom(X,-)).
$$
Dually, one considers proper $\catx$-coresolutions
to define the functors $\axext^n(-,M)$.
Consult Section~\ref{sec1} for  precise definitions.

In this article we investigate relative cohomology theories that 
arise from dualities with respect to semidualizing modules:
when $R$ is commutative and noetherian, a finitely generated $R$-module $C$ 
is \emph{semidualizing} when $\aext^{\geq 1}(C,C)=0$ and
$\ahom(C,C)\cong R$.
Examples include projective $R$-modules of rank 1 and,
when $R$ is a Cohen-Macaulay ring of finite Krull dimension that is a homomorphic image of a
Gorenstein ring, a dualizing module.

A semidualizing $R$-module $C$ gives rise to several 
distinguished 
classes of modules.  For instance, one has the 
class $\catpc$ of $C$-projective modules and the class $\catgpc$
of $\text{G}_C$-projective modules,  which we use for resolutions.
For coresolutions, we consider the class $\catic$ of $C$-injective modules  and the class
$\catgic$ of $\text{G}_C$-injective modules.
As there is no risk of confusion in these cases, the corresponding relative
cohomology functors are denoted 
$\pcext^n(-,-)$, 
$\gpcext^n(-,-)$, 
$\icext^n(-,-)$ and
$\gicext^n(-,-)$.
Detailed definitions can be found in Section~\ref{sec8}.

Our investigation into these functors focuses on two questions:  
What conditions on a pair of modules $(M,N)$ guarantee that
the corresponding outputs of two of these functors are isomorphic?
And when are these functors different?

As to the first question,
Section~\ref{sec7} focuses on the issue of balance,
motivated by the fact that one can compute the 
``absolute'' cohomology $\aext^n(M,N)$ in terms of
a projective resolution of $M$ or an injective
resolution of $N$.  This section begins with Example~\ref{ctrex2}
which shows
that the naive version of balance for relative cohomology fails in general:
even if $\pcpd_R(M)$ and $\icid_R(N)$ are both finite, one can have
$\pcext^n(M,N)\ncong\icext^n(M,N)$ and 
$\gpcext^n(M,N)\ncong\gicext^n(M,N)$.  
It turns out that the correct balance result in this setting,
stated next, uses
coresolutions with respect to the semidualizing module 
$\da{C}=\ahom(C,D)$ where $D$ is a dualizing module.
This result is contained in Proposition~\ref{balance02b} and
Theorem~\ref{balance02a}.

\begin{intthm} \label{introbalance}
Let $R$ be a Cohen-Macaulay ring with a dualizing module,
and let $C$, $M$ and $N$ be $R$-modules with $C$ semidualizing.
\begin{enumerate}[\quad\rm(a)]
\item \label{introbalanceitem1}
If $\pcpd_R(M)<\infty$ and $\icdid_R(N)<\infty$, then there is an isomorphism
$$\pcext^n(M,N)\cong\icdext^n(M,N)$$
for each integer $n$.
\item \label{introbalanceitem2}
If $\gpcpd_R(M)<\infty$ and $\gicdid_R(N)<\infty$, then there is an isomorphism
$$\gpcext^n(M,N)\cong\gicdext^n(M,N)$$
for each integer $n$.
\end{enumerate}
\end{intthm}

In addition,
Section~\ref{sec5} gives conditions that yield isomorphisms
$\pcext^n(M,N)\cong\gpcext^n(M,N)$ and 
$\icext^n(M,N)\cong\gicext^n(M,N)$.  See Propositions~\ref{detect01}
and~\ref{detect01a}.  

Section~\ref{sec08} deals with the 
even more interesting question of the differences between these 
functors.  The next result summarizes our findings from this section and shows that
each reasonably comparable pair of relative cohomology functors is
distinct.

\begin{intthm} \label{intronotiso}
Let $(R,\m)$ be a local ring,
and let $B,C$ be semidualizing $R$-modules. 
\begin{enumerate}[\quad\rm(a)]
\item \label{intronotisoitem1}
Assume $C\ncong R$.  
Then one has 
\begin{gather*}
\pcext(-,-)\ncong\aext(-,-)\ncong\gpcext(-,-)\\
\icext(-,-)\ncong\aext(-,-)\ncong\gicext(-,-).
\end{gather*}
If there exist $y,z\in\m$ such that $\ann_R(y)=zR$ and
$\ann_R(z)=yR$, then 
\begin{gather*}
\pcext(-,-)\ncong\gpcext(-,-) \qquad
\icext(-,-)\ncong\gicext(-,-).
\end{gather*}
\item \label{intronotisoitem2}
Assume $\gpcpd_R(B)<\infty$ and $C\ncong B$.  
Then one has 
\begin{gather*}
\pcext(-,-)\ncong\pbext(-,-)\ncong\gpcext(-,-) \\
\gibext(-,-)\ncong\icext(-,-)\ncong\ibext(-,-)\ncong\gicext(-,-).
\end{gather*}
If $\depth(R)\geq 1$, then 
\begin{gather*}
\gpcext(-,-)\ncong\gpbext(-,-) \qquad
\gicext(-,-)\ncong\gibext(-,-).
\end{gather*}
If $C$ admits a proper $\catgpb$-resolution, then 
$$\pcext(-,-)\ncong\gpbext(-,-).$$
\end{enumerate}
\end{intthm}

As an aid for some of the computations in Theorem~\ref{intronotiso}
we utilize a Yoneda-type characterization of relative cohomology modules.
This is the subject of Section~\ref{sec2}.  In particular, the following result
is contained in Theorem~\ref{extincl}.

\begin{intthm} \label{introextincl}
Let $M$ and $N$ be $R$-modules.
\begin{enumerate}[\quad\rm(a)]
\item \label{introextinclitem1}
If $M$ admits a proper $\catx$-resolution, then 
$\xaext(M,N)$ is in bijection with the 
set of equivalence classes of  sequences
$0\to N\to T\to M\to 0$ that are exact and
$\hom_R(\catx,-)$-exact.
\item \label{introextinclitem2}
If $N$ admits a proper $\caty$-coresolution, then 
$\ayext(M,N)$ is in bijection with the 
set of equivalence classes of  sequences
$0\to N\to T\to M\to 0$ that are exact and
$\hom_R(-,\caty)$-exact.
\end{enumerate}
\end{intthm}

\section{Categories, Resolutions, and Relative Cohomology}\label{sec1}

We begin with some notation and terminology for use throughout this paper.

\begin{notationdefinition} \label{notation01}
Throughout this work
$R$ is a commutative ring.
Write $\catm=\catm(R)$ 
for the category of 
$R$-modules,
and write $\catp=\catp(R)$, $\catf=\catf(R)$ and $\cati=\cati(R)$ 
for the subcategories of projective, flat and injective
$R$-modules, respectively.
We use the term ``subcategory'' to mean a ``full, additive, and essential (closed under isomorphisms) 
subcategory.''
If $\catx$ is a subcategory of $\catm$, then $\catx^f$ is the  subcategory
of finitely generated modules in $\catx$.  
\end{notationdefinition}

\begin{defn} \label{notation01a}
We fix subcategories $\catx$, $\caty$, $\catw$, and $\catv$  of $\catm$ such that
$\catw\subseteq\catx$
and $\catv \subseteq\caty$.
Write $\catx\perp\caty$ 
if $\ext_R^{\geq1}(X,Y)=0$ for each module $X$ in $\catx$ and each module $Y$ in $\caty$.
For a module $M$ in $\catm$, write $M\perp\caty$ (resp., $\catx\perp M$)
if $\ext_R^{\geq1}(M,Y)=0$ for each module  $Y$ in $\caty$
(resp., if $\ext_R^{\geq1}(X,M)=0$ for each module  $X$ in $\catx$).
We say that $\catw$ is a \emph{cogenerator} for $\catx$ if,
for each module $X$ in $\catx$, there exists an exact sequence 
$0\to X\to W\to X'\to 0$
such that $W$ is in $\catw$ and $X'$ is in $\catx$.
The subcategory
$\catw$ is an \emph{injective cogenerator} for $\catx$ if 
$\catw$ is a cogenerator for $\catx$ and $\catx\perp\catw$.
The terms \emph{generator} and \emph{projective generator} 
are defined dually.
\end{defn}

\begin{defn} \label{notation07}
An \emph{$R$-complex} is a sequence of 
$R$-module homomorphisms 
$$X =\cdots\xra{\partial^X_{n+1}}X_n\xra{\partial^X_n}
X_{n-1}\xra{\partial^X_{n-1}}\cdots$$
such that $ \partial^X_{n-1}\partial^X_{n}=0$ for each integer $n$; the
$n$th \emph{homology module} of $X$ is
$\HH_n(X)=\Ker(\partial^X_{n})/\im(\partial^X_{n+1})$.
We frequently identify $R$-modules with complexes concentrated in degree 0.
The  \emph{suspension} (or \emph{shift}) of
$X$, denoted $\shift X$, is the complex with
$(\shift X)_n=X_{n-1}$ and $\partial_n^{\shift X}=-\partial_{n-1}^X$.

The complex $X$ is \emph{$\hom_R(\catx,-)$-exact} if the complex
$\hom_R(X',X)$ is exact for each module $X'$ in $\catx$.  
Dually, it is \emph{$\hom_R(-,\catx)$-exact} if the complex
$\hom_R(X,X')$ is exact for each module $X'$ in $\catx$. 
It is \emph{$-\otimes_R\catx$-exact} if the complex
$X'\otimes_R X$ is exact for each module $X'$ in $\catx$. 
\end{defn}

\begin{defn} \label{notation07a}
Let $X,Y$ be $R$-complexes.
The  Hom-complex $\hom_R(X,Y)$ is the $R$-complex defined as
$\hom_R(X,Y)_n=\prod_p\hom_R(X_p,Y_{p+n})$
with $n$th differential $\partial_n^{\hom_R(X,Y)}$ given by
$\{f_p\}\mapsto \{\partial^{Y}_{p+n}f_p-(-1)^nf_{n-1}\partial^X_p\}$.
A \emph{morphism} 
is an element of $\ker(\partial_0^{\hom_R(X,Y)})$.
Two morphisms $\alpha,\alpha'\colon X\to Y$
are \emph{homotopic}, written $\alpha\sim\alpha'$, if the difference
$\alpha-\alpha'$ is in $\im(\partial_1^{\hom_R(X,Y)})$.
The morphism $\alpha$ is a
\emph{homotopy equivalence} if there is a morphism
$\beta\colon Y\to X$ such that 
$\beta \alpha \sim\id_{X}$ and
$\alpha\beta\sim\id_{Y}$.

A morphism of complexes $\alpha\colon X\to Y$
induces homomorphisms on homology modules
$\HH_n(\alpha)\colon\HH_n(X)\to\HH_n(Y)$, and $\alpha$ is a
\emph{quasiisomorphism} when each $\HH_n(\alpha)$ is bijective.
The \emph{mapping cone} of $\alpha$ is the complex
$\cone(\alpha)$ defined as
$\cone(\alpha)_n=Y_n\oplus X_{n-1}$
with $n$th differential 
$\partial^{\cone(\alpha)}_n 
= \Bigl(\begin{smallmatrix}\partial_{n}^{Y} & \alpha_{n-1} \\ 0 & -\partial_{n-1}^{X}
\end{smallmatrix} \Bigr)$.
Recall that $\alpha$ is a quasiisomorphism if and only if $\cone(\alpha)$ is
exact.
\end{defn}

\begin{defn} \label{notation03}
An $R$-complex $X$ is \emph{bounded} if $X_n=0$ for $|n|\gg 0$. When
$X_{-n}=0=\HH_n(X)$ for all $n>0$, the natural map
$X\to\HH_0(X)\cong M$ is a quasiisomorphism.  In this event, $X$ is an
\emph{$\catx$-resolution} of $M$ if each $X_n$ is in $\catx$, and
the exact sequence
$$X^+ = \cdots\xra{\partial^X_{2}}X_1
\xra{\partial^X_{1}}X_0\to M\to 0$$ is the \emph{augmented
$\catx$-resolution} of $M$ associated to $X$.
We write ``projective resolution'' in lieu of
``$\catp$-resolution''.
The \emph{$\catx$-projective dimension} of $M$ is the quantity
$$\xpd(M)=\inf\{\sup\{n\geq 0\mid X_n\neq 0\}\mid \text{$X$ is an
$\catx$-resolution of $M$}\}.$$ 
The modules of $\catx$-projective dimension 0 are
the nonzero modules in $\catx$.
We let
$\finrescatx$ denote
the subcategory of $R$-modules $M$  with $\xpd(M)<\infty$.
One checks easily that $\finrescatx$ is additive and contains $\catx$.

The terms \emph{$\caty$-coresolution} and \emph{$\caty$-injective dimension}
are defined dually.  The \emph{augmented 
$\caty$-coresolution} associated to a $\caty$-coresolution $Y$ is denoted $^+Y$,
and the $\caty$-injective dimension of $M$ is denoted $\yid(M)$.
The subcategory of $R$-modules $N$  with $\yid(N)<\infty$ is denoted
$\fincorescaty$; it is additive and contains $\caty$.

Following much of the literature, we write ``injective resolution'' in lieu of
``$\cati$-coresolution'' and set $\pd=\catp\text{-}\pd$ and $\id=\cati\text{-}\id$.
\end{defn}

\begin{defn} \label{notation05}
An $\catx$-resolution $X$ is \emph{proper} if
the augmented resolution $X^+$ is $\hom_R(\catx,-)$-exact. We let
$\proprescatx$ denote
the subcategory of $R$-modules admitting a proper 
$\catx$-resolution.
One checks readily that $\proprescatx$ is additive and contains
$\catx$.
\emph{Proper coresolutions} are defined dually.
The subcategory of $R$-modules admitting a proper 
$\caty$-coresolution is denoted $\propcorescaty$; 
it is additive and contains $\caty$.
\end{defn}

The next lemmata are standard or have standard proofs:
for~\ref{perp03} see~\cite[pf.~of (2.3)]{auslander:htmcma};
for~\ref{gencat01} see~\cite[pf.~of (2.1)]{auslander:htmcma};
for~\ref{rel01} argue as in~\cite[(4.3)]{avramov:aratc} and~\cite[(1.8)]{holm:ghd}; and
for the ``Horseshoe Lemma'' \ref{horseshoe01} 
see~\cite[(4.5)]{avramov:aratc}
and~\cite[pf.~of (8.2.1)]{enochs:rha}.

\begin{lem} \label{perp03}
Let $0\to M_1\to M_2\to M_3\to 0$ be an exact sequence of $R$-modules.
\begin{enumerate}[\quad\rm(a)]
\item \label{perp03item1}
If $M_3\perp\catx$, then $M_1\perp\catx$ if and only if $M_2\perp\catx$.
If $M_1\perp\catx$  and  $M_2\perp\catx$, 
then $M_3\perp\catx$
if and only if the given sequence is  $\ahom(-,\catx)$-exact.
\item \label{perp03item2}
If $\catx\perp M_1$, then $\catx\perp M_2$ if and only if $\catx\perp M_3$.
If $\catx\perp M_2$  and  $\catx\perp M_3$, 
then $\catx\perp M_1$
if and only if the given sequence is  $\ahom(\catx,-)$-exact. \qed
\end{enumerate}
\end{lem}

\begin{lem} \label{gencat01}
If $\catx\perp\caty$, then  $\catx\perp\finrescaty$ and $\fincorescatx\perp\caty$.
\qed
\end{lem}

\begin{lem} \label{rel01}
Let $M,M',N,N'$ be $R$-modules.
\begin{enumerate}[\quad\rm(a)]
\item \label{rel01item1}
Let $P\xra{\rho} M$ be a projective resolution.
Assume that $M$ admits a proper $\catw$-resolution
$W\xra{\gamma} M$ and $M'$ admits a proper $\catx$-resolution
$X'\xra{\gamma'} M'$.
For each homomorphism $f\colon M\to M'$ there exist
morphisms $\ol{f}\colon W\to X'$
and $\wti{f}\colon P\to X'$ unique up to homotopy
such that $f\gamma=\gamma'\ol{f}$ and $f\rho=\gamma'\wti{f}$.
If $f$ is an isomorphism, then $\ol{f}$ and $\wti{f}$ are quasiisomorphisms.
If $f$ is an isomorphism and $\catx=\catw$, then  $\ol{f}$ is a
homotopy equivalence.
\item \label{rel01item2}
Let $N'\xra{\phi'} I'$ be an injective resolution.
Assume that $N$  admits a proper $\caty$-coresolution
$N\xra{\delta} Y$ and $N'$  admits a proper $\catv$-coresolution
$N'\xra{\delta'} V'$.
For each homomorphism $g\colon N\to N'$ there exist
morphisms $\ol{g}\colon Y\to V'$
and $\wti{g}\colon Y\to I'$ unique up to homotopy
such that $\ol{g}\delta=\delta' g$ and $\wti{g}\delta=\phi' g$.
If $g$ is an isomorphism, then $\ol{g}$ and $\wti{g}$ are quasiisomorphisms.
If $g$ is an isomorphism and $\catv=\caty$, then $\ol{g}$  is a
homotopy equivalence. \qed
\end{enumerate}
\end{lem}

\begin{lem} \label{horseshoe01}
Let $0\to M'\to M\to M''\to 0$ be an exact sequence of $R$-modules.
\begin{enumerate}[\quad\rm(a)]
\item \label{horseshoe01item1}
Assume that $M'$ and $M''$ 
admit proper $\catx$-resolutions $X'\xra{\simeq} M'$ and $X''\xra{\simeq}M''$
and that the given sequence
is $\ahom(\catx,-)$-exact. Then $M$
admits a proper $\catx$-resolution $X\xra{\simeq} M$
such that there exists a commutative diagram 
$$
\xymatrix{
0\ar[r] & 
X' \ar[r]\ar[d]_{\simeq} 
& X \ar[r]\ar[d]_{\simeq} 
& X'' \ar[r]\ar[d]_{\simeq} & 0 \\
0\ar[r] & M' \ar[r] & M \ar[r] & M'' \ar[r] & 0. 
}
$$
whose top row is  degreewise
split exact.  
\item \label{horseshoe01item2}
Assume that $M'$ and $M''$ 
admit proper $\caty$-coresolutions $M'\xra{\simeq} Y'$ and $M''\xra{\simeq}Y''$
and that the given sequence
is $\ahom(-,\caty)$-exact. Then $M$ admits a
proper $\caty$-coresolution $M\xra{\simeq} Y$
such that there exists a commutative diagram 
$$
\xymatrix{
0\ar[r] & M' \ar[r]\ar[d]_{\simeq} & M \ar[r]\ar[d]_{\simeq} & M'' \ar[r]\ar[d]_{\simeq} & 0 \\
0\ar[r] 
& Y' \ar[r]
& Y \ar[r]
& Y'' \ar[r] & 0 
}
$$
whose bottom row is  degreewise
split exact.
\qed
\end{enumerate}
\end{lem}

\begin{defn} \label{rel02}
Let $M,M',N,N'$ be $R$-modules
equipped with homomorphisms
$f\colon M\to M'$ and $g\colon N\to N'$.
Assume that $M$ admits a proper $\catx$-resolution 
$X\xra{\gamma} M$, and define
the \emph{$n$th relative $\catx R$-cohomology module} as
$$\xaext^n(M,N)=\HH_{-n}(\ahom(X,N))$$
for each integer $n$.
If $M'$ also admits a proper $\catx$-resolution 
$X'\xra{\gamma'} M'$, then 
let $\ol{f}\colon X\to X'$ be a chain map such that $f\gamma=\gamma'\ol{f}$
as in Lemma~\ref{rel01} and define
\begin{gather*}
\xaext^n(f,N)=\HH_{-n}(\ahom(\ol{f},N))\colon\xaext^n(M',N)\to\xaext^n(M,N) \\
\xaext^n(M,g)=\HH_{-n}(\ahom(X,g))\colon\xaext^n(M,N)\to\xaext^n(M,N'). 
\end{gather*}
The \emph{$n$th relative $R\caty$-cohomology}
$\ayext^n(-,-)$ is defined dually.
\end{defn}

\begin{disc} \label{rel03}
Lemma~\ref{rel01} shows that Definition~\ref{rel02} yields well-defined bifunctors
\begin{align*}
\xaext^n(-,-)&\colon\proprescatx\times\catm\to \catm & \text{and} & &
\ayext^n(-,-)&\colon\catm\times \propcorescaty\to \catm
\end{align*}
and one checks the following natural equivalences readily.
\begin{gather*}
\xaext^{\geq 1}(\catx,-)=0=\ayext^{\geq 1}(-,\caty) \\
\xaext^0(-,-)\cong\ahom(-,-)|_{\proprescatx\times\catm} \\
\ayext^0(-,-)\cong\ahom(-,-)|_{\catm\times \propcorescaty} \\
\ext^n_{\catp\catm}(-,-)\cong\aext^n(-,-)\cong\ext^n_{\catm\cati}(-,-)
\end{gather*}
\end{disc}

\begin{defn} \label{rel02a}
Let $M,N$ be $R$-modules.
Let $P\xra{\rho}M$ be a projective resolution.
Assume that $M$ admits a proper $\catw$-resolution 
$W\xra{\gamma} M$ and a proper $\catx$-resolution 
$X\xra{\gamma'} M$. Let $\ol{\id_M}\colon W\to X$
and $\wti{\id_M}\colon P\to X$ be  quasiisomorphisms
such that $\gamma=\gamma'\ol{\id_M}$ and $\rho=\gamma'\wti{\id_M}$,
as in 
Lemma~\ref{rel01}\eqref{rel01item1}, and set
\begin{align*}
\xwacomp^n(M,N)=\HH_{-n}(\ahom(\ol{\id_M},N))
&\colon\xaext^n(M,N)\to\waext^n(M,N) \\
\xrcomp^n(M,N)=\HH_{-n}(\ahom(\wti{\id_M},N))
&\colon\xaext^n(M,N)\to\ext_R^n(M,N).
\end{align*}
On the other hand, if $N$ admits a proper $\caty$-coresolution 
and a proper $\catv$-coresolution, then the following maps are defined dually 
\begin{align*}
\ayvcomp^n(M,N)&\colon\ayext^n(M,N)\to\avext^n(M,N) \\
\rycomp^n(M,N)&\colon\ayext^n(M,N)\to\aext^n(M,N). 
\end{align*}
\end{defn}

\begin{disc}
Lemma~\ref{rel01} shows that Definition~\ref{rel02a} describes
well-defined natural transformations
\begin{align*}
\xwacomp^n(-,-)&\colon\xaext^n(-,-)|_{(\proprescatw\cap \proprescatx)\times\catm}\to 
\waext^n(-,-)|_{(\proprescatw\cap \proprescatx)\times\catm} \\
\xrcomp^n(-,-)&\colon\xaext^n(-,-)\to 
\aext^n(-,-)|_{\proprescatx\times\catm} \\
\ayvcomp^n(-,-)&\colon\ayext^n(-,-)|_{\catm \times(\propcorescatv\cap \propcorescaty)}\to 
\avext^n(-,-)|_{\catm \times(\propcorescatv\cap \propcorescaty)} \\
\rycomp^n(-,-)&\colon\ayext^n(-,-)\to 
\aext^n(-,-)|_{\catm \times\propcorescaty}
\end{align*}
independent of resolutions and liftings.
Note that the left-exactness of $\hom_R(-,-)$ implies that
 each of these transformations is a natural isomorphism when $n\leq 0$.
\end{disc}

Lemma~\ref{horseshoe01} yields the following long exact sequences 
as in~\cite[(4.4),(4.6)]{avramov:aratc}.

\begin{lem} \label{notation06a}
Let $M$ and $N$ be $R$-modules, and consider an exact  sequence 
$$\mathbf{L}=\quad 0\to L'\xra{f'} L\xra{f} L''\to 0.$$
\begin{enumerate}[\quad\rm(a)]
\item \label{06aitem1}
Assume that the sequence $\mathbf{L}$ is $\ahom(\catx,-)$-exact.
If $M$ is in $\proprescatx$, then $\mathbf{L}$ induces a functorial long exact sequence
\begin{align*}
\cdots \to
& \xaext^n(M,L') \xra{\xaext^n(M,f')}
 \xaext^n(M,L) \xra{\xaext^n(M,f)}
  \\
& \xaext^n(M,L'')\xra{\eth^n_{\catx\!\cata}(M,\mathbf{L})}
 \xaext^{n+1}(M,L') \xra{\xaext^{n+1}(M,f')}
 \cdots. 
\end{align*}
\item \label{06aitem2}
Assume that the sequence $\mathbf{L}$ is $\ahom(\catx,-)$-exact.
If the modules $L',L,L''$ are in 
$\proprescatx$, then $\mathbf{L}$ induces a functorial long exact sequence
\begin{align*}
\cdots \to
& \xaext^n(L'',N) \xra{\xaext^n(f,N)}
 \xaext^n(L,N) \xra{\xaext^n(f',N)}
  \\
& \xaext^n(L',N)\xra{\eth^n_{\catx\!\cata}(\mathbf{L},N)}
 \xaext^{n+1}(L'',N) \xra{\xaext^{n+1}(f,N)}
 \cdots. 
\end{align*}
\item \label{06aitem3}
Assume that the sequence $\mathbf{L}$ is  $\ahom(-,\caty)$-exact.
If $N$ is in $\propcorescaty$, then $\mathbf{L}$ induces a functorial long exact sequence
\begin{align*}
\cdots \to
& \ayext^n(L'',N) \xra{\ayext^n(f,N)}
 \ayext^n(L,N) \xra{\ayext^n(f',N)}
  \\
& \ayext^n(L',N)\xra{\eth^n_{\cata\caty}(\mathbf{L},N)}
 \ayext^{n+1}(L'',N) \xra{\ayext^{n+1}(f,N)}
 \cdots. 
\end{align*}
\item \label{06aitem4}
Assume that the sequence $\mathbf{L}$ is  $\ahom(-,\caty)$-exact.
If the modules $L',L,L''$ are in 
$\propcorescaty$, then $\mathbf{L}$ induces a functorial long exact sequence
\begin{align*}
&&&&&&&&\cdots \to
& \ayext^n(M,L') \xra{\ayext^n(M,f')}
 \ayext^n(M,L) \xra{\ayext^n(M,f)}
  \\
&&&&&&&&& \ayext^n(M,L'')\xra{\eth^n_{\cata\caty}(M,\mathbf{L})}
 \ayext^{n+1}(M,A') \xra{\ayext^{n+1}(M,f')}
 \cdots. &&&& \qed
\end{align*}
\end{enumerate}
\end{lem}

\section{Relative Cohomology and Extensions} \label{sec2}

In this section, we compare relative cohomology modules with 
sets of equivalence classes of module extensions, as in the classical
Yoneda setting.  

\begin{notationdefinition} \label{yoneda}
Let $N$ and $M$ be $R$-modules.
An \emph{extension} of $M$ by $N$ is an exact sequence
$$0\to N\to T\to M\to 0$$
and this is \emph{equivalent} to a second extension
$0\to N\to T'\to M\to 0$ if there exists a homomorphism
$\tau\colon T\to T'$ making the following diagram commute
$$
\xymatrix{
0 \ar[r] & N \ar[r] \ar[d]_{\id_N} & T  \ar[r] \ar[d]_{\tau} & M \ar[r] \ar[d]_{\id_M} & 0 \\
0 \ar[r] & N \ar[r] & T'  \ar[r]  & M \ar[r]  & 0 
}$$
We set 
\begin{align*}
\eext_R(M,N)&=\{\text{equivalence classes of extensions
of $M$ by $N$}\} \\
\xeext(M,N)&=\{\text{equivalence classes of 
$\hom_R(\catx,-)$-exact extensions
of $M$ by $N$}\} \\
\yeext(M,N)&=\{\text{equivalence classes of 
$\hom_R(-,\caty)$-exact extensions
of $M$ by $N$}\} 
\end{align*}
\end{notationdefinition}

\begin{disc}  \label{yoneda'}
Because of the containments $\catw\subseteq\catx$ and $\catv\subseteq\caty$
there are inclusions
$\weext(M,N)\subseteq\xeext(M,N)\subseteq\eext(M,N)$
and $\veext(M,N)\subseteq\yeext(M,N)\subseteq\eext(M,N)$.

There exists a bijection 
$\xi_{RMN}\colon \ext^1_R(M,N)\to\eext_R(M,N)$ 
whose construction we recall from~\cite[Ch.~7]{rotman:iha}. Let $P\xra{\simeq} M$ be a projective
resolution.  Each element in $\ext_R^1(M,N)$ is represented by 
a homomorphism $\alpha\colon P_1\to N$ such that $\alpha\partial_2^P=0$, and
each such $\alpha$ induces a map $\overline{\alpha}\colon P_1/\im(\partial^P_2)\to N$.
Taking a pushout yields the following commutative diagram with exact rows
$$
\xymatrix{
0 \ar[r] 
& P_1/\im(\partial^P_2) \ar[r]^-{\overline{\partial^P_1}} \ar[d]_{\ol{\alpha}} \ar@{}[rd]|>>{\lrcorner} 
& P_0  \ar[r] \ar[d]_{\tau} & M \ar[r] \ar[d]_{\id_M} & 0 \\
0 \ar[r] & N \ar[r] & T \ar[r]  & M \ar[r]  & 0 
}$$
and $\xi_{RMN}([\alpha])$ is the equivalence class of the bottom row
of this diagram. 

Dually, one constructs a bijection 
$\xi'_{RMN}\colon \ext^1_R(M,N)\to\eext_R(M,N)$ 
using an injective resolution of $N$.
\end{disc}

The next result extends the construction of Remark~\ref{yoneda'} to the relative setting
and contains Theorem~\ref{introextincl} from the introduction.
The  connecting maps 
$\xrcomp^1$, $\xwacomp^1$, $\rycomp^1$ and $\ayvcomp^1$
are described
in Definition~\ref{rel02a}.

\begin{thm} \label{extincl}
Let $M$ and $N$ be $R$-modules.
\begin{enumerate}[\quad\rm(a)]
\item \label{extinclitem1}
Assume that $M$ admits a proper $\catx$-resolution.
There is then a bijection
$\xi_{\catx MN}\colon \xaext^1(M,N)\to\xeext(M,N)$ 
making the following diagram commute
$$\xymatrix{
\xaext^1(M,N) \ar[r]^-{\xi_{\catx MN}} \ar[d]_{\xrcomp^1(M,N)} & \xeext(M,N) \ar[d] \\
\aext^1(M,N) \ar[r]^-{\xi_{RMN}} & \eext(M,N)
}$$
where the rightmost vertical arrow is the natural inclusion.
In particular, the comparison 
map $\xrcomp^1(M,N)\colon\xaext^1(M,N)\to\aext^1(M,N)$ is
injective.
\item \label{extinclitem1'}
Assume that $M$ admits a proper $\catx$-resolution
and a proper $\catw$-resolution. The following diagram commutes
where the rightmost vertical arrow is the natural inclusion
$$\xymatrix{
\xaext^1(M,N) \ar[r]^-{\xi_{\catx MN}} \ar[d]_{\xwacomp^1(M,N)} & \xeext(M,N) \ar[d] \\
\waext^1(M,N) \ar[r]^-{\xi_{\catw MN}} & \weext(M,N).
}$$
In particular the comparison map $\xwacomp^1(M,N)$ is injective.
\item \label{extinclitem2}
Assume that $N$ admits a proper $\caty$-coresolution.
There is then a bijection
$\xi'_{\caty MN}\colon \ayext^1(M,N)\to\yeext(M,N)$ 
making the following diagram commute
$$\xymatrix{
\ayext^1(M,N) \ar[r]^-{\xi'_{\caty MN}} \ar[d]_{\rycomp^1(M,N)} & \yeext(M,N) \ar[d] \\
\aext^1(M,N) \ar[r]^-{\xi'_{RMN}} & \eext(M,N)
}$$
where the rightmost vertical arrow is the natural inclusion.
In particular, the comparison 
map $\rycomp^1(M,N)\colon\ayext^1(M,N)\to\aext^1(M,N)$ is
injective.
\item \label{extinclitem2'}
Assume that $N$ admits a proper $\caty$-coresolution
and a proper $\catv$-coresolution. The following diagram commutes
where the rightmost vertical arrow is the natural inclusion
$$\xymatrix{
\ayext^1(M,N) \ar[r]^-{\xi'_{\caty MN}} \ar[d]_{\ayvcomp^1(M,N)} & \yeext(M,N) \ar[d] \\
\avext^1(M,N) \ar[r]^-{\xi'_{\catv MN}} & \veext(M,N).
}$$
In particular the comparison map $\ayvcomp^1(M,N)$ is injective.
\end{enumerate}
\end{thm}

\begin{proof}
Our proof is modeled
on the arguments of~\cite[Ch.~7]{rotman:iha}; instead of rewriting much of the 
work there, we simply sketch the proof, indicating 
how  $\hom_R(\catx,-)$-exactness is detected and used.

\eqref{extinclitem1} Let $X\xra{\simeq} M$ be a proper
$\catx$-resolution and set $\ol{X}=X_1/\im(\partial^X_2)$.
Each element $[\alpha]\in\xaext^1(M,N)$ is represented by 
a homomorphism $\alpha\colon X_1\to N$ such that $\alpha\partial_2^X=0$,
and each such $\alpha$ induces a map $\overline{\alpha}\colon \ol{X}\to N$.
Taking a pushout yields the following commutative diagram with exact rows
\begin{equation} \label{yoneda01}
\begin{split}
\xymatrix{
& 0 \ar[r] 
& \ol{X} \ar[r]^-{\overline{\partial^X_1}} \ar[d]_{\ol{\alpha}} \ar@{}[rd]|>>{\lrcorner} 
& X_0  \ar[r] \ar[d]_{\tau} & M \ar[r] \ar[d]_{\id_M} & 0 \\
\zeta= & 0 \ar[r] & N \ar[r] & T \ar[r]^-{\pi}  & M \ar[r]  & 0. 
}
\end{split}
\end{equation}
We claim that the bottom row of this diagram is $\hom_R(\catx,-)$-exact.
To see this, first note that the properness of the resolution $X\xra{\simeq} M$
implies that the top row of~\eqref{yoneda01} is $\hom_R(\catx,-)$-exact.
Fix an $R$-module $X'$ in $\catx$ and apply $\hom_R(X',-)$ to the
diagram~\eqref{yoneda01} to yield the next commutative diagram
with exact rows
$$
\xymatrix{
0 \ar[r] 
& \hom_R(X',\ol{X}) \ar[r] \ar[d] 
& \hom_R(X',X_0)  \ar[rr] \ar[d] && \hom_R(X',M) \ar[r] \ar[d]^{\id_{\hom_R(X',M)}} & 0 \\
0 \ar[r] & \hom_R(X',N) \ar[r] & \hom_R(X',T) \ar[rr]^-{\hom_R(X',\pi)}  && \hom_R(X',M) . 
}
$$
An easy diagram chase shows that the map $\hom_R(X',\pi)$ is surjective, as desired.

We define $\xi_{\catx MN}([\alpha])$ to be the equivalence class $[\zeta]$
of the bottom row
of the diagram~\eqref{yoneda01}. 
One now verifies readily
(as in the proof of  the classical result in~\cite[Ch.~7]{rotman:iha})
that this yields a well-defined function $\xaext^1(M,N)\to\xeext(M,N)$. 

To show that $\xi_{\catx MN}$ is bijective, we construct an inverse.  Fix
a $\ahom(\catx,-)$-exact sequence 
$\zeta=(0\to N\to T\to M\to 0)$.  
A standard lifting procedure as in Lemma~\ref{rel01}\eqref{rel01item1}
yields the next commutative diagram with exact rows
$$
\xymatrix{
X_2 \ar[r]^-{\partial^X_2} \ar[d]
& X_1 \ar[r]^-{\partial^X_1} \ar[d]_{\alpha}
& X_0  \ar[r] \ar[d] & M \ar[r] \ar[d]_{\id_M} & 0 \\
0 \ar[r] & N \ar[r] & T \ar[r]  & M \ar[r]  & 0. 
}
$$
The map $\alpha$ is thus a degree-1 cycle in $\hom_R(X,N)$ and so gives rise to a
cohomology class $[\alpha]\in \xaext^1(M,N)$.
Again, one verifies that the assignment $[\zeta]\mapsto[\alpha]$
describes a well-defined function $\xeext(M,N)\to \xaext^1(M,N)$,
and that this function is a two-sided inverse for $\xi_{\catx MN}$;
the reader may find~\cite[(7.18)]{rotman:iha} to be helpful.

The proof of part~\eqref{extinclitem1} will be compete once we 
verify $\xi_{RMN}\xrcomp^1=\xi_{\catx MN}$.
Let $P\xra{\simeq} M$ be a projective
resolution and set $\ol{P}=P_1/\im(\partial^P_2)$. 
Lemma~\ref{rel01}\eqref{rel01item1} yields 
the next commutative diagram
with exact rows
$$
\xymatrix{
\cdots \ar[r]^-{\partial^P_3} & P_2 \ar[r]^-{\partial^P_2} \ar[d]_{f_2}
& P_1 \ar[r]^-{\partial^P_1} \ar[d]_{f_1}
& P_0  \ar[r] \ar[d]_{f_0} 
& M \ar[r] \ar[d]_{\id_M} & 0 \\
\cdots \ar[r]^-{\partial^X_3} & X_2 \ar[r]^-{\partial^X_2} 
& X_1 \ar[r]^-{\partial^X_1}
& X_0  \ar[r]  & M \ar[r]  & 0 
}
$$
which in turn induces another commutative diagram
with exact rows
\begin{equation} \label{yoneda03}
\begin{split}
\xymatrix{
0 \ar[r]
& \ol{P} \ar[r]^-{\ol{\partial^P_1}} \ar[d]_{\ol{f_1}}
& P_0  \ar[r] \ar[d]_{f_0} 
& M \ar[r] \ar[d]_{\id_M} 
& 0 \\
0 \ar[r]
& \ol{X} \ar[r]^-{\ol{\partial^X_1}}
& X_0  \ar[r]  
& M \ar[r]  
& 0 
}
\end{split}
\end{equation}
Given $[\alpha]\in\xaext^1(M,N)$, construct the extension
$\zeta$ as above.
The diagram~\eqref{yoneda03} combines with~\eqref{yoneda01} to yield the next
diagram
\begin{equation} \label{yoneda02}
\begin{split}
\xymatrix{
& 0 \ar[r] 
& \ol{P} \ar[r]^-{\overline{\partial^P_1}} \ar[d]_{\ol{\alpha}\ol{f_1}} 
& P_0  \ar[r] \ar[d]_{\tau f_0} & M \ar[r] \ar[d]_{\id_M} & 0 \\
\zeta= & 0 \ar[r] & N \ar[r] & T \ar[r]^-{\pi}  & M \ar[r]  & 0. 
}
\end{split}
\end{equation}
It follows from Definition~\ref{rel02a} that 
$\xrcomp^1([\alpha])=[\alpha f_1]$.
From~\cite[(7.18)]{rotman:iha} one concludes
$\xi_{RMN}([\alpha f_1])=[\zeta]$,
and this yields the first equality in the following sequence
$$
\xi_{RMN}(\xrcomp^1([\alpha]))
=[\zeta]
=\xi_{\catx MN}([\alpha])
$$
while the second equality is by definition.  This
completes the proof of part~\eqref{extinclitem1}.

Part~\eqref{extinclitem1'} is proved as in the previous paragraph,
using a proper $\catw$-resolution in place of the projective resolution $P\to M$.
The proofs
of~\eqref{extinclitem2} and~\eqref{extinclitem2'} are dual.  
\end{proof}

\section{Categories of Interest} \label{sec8}

In this section we discuss the categories whose relative cohomology theories
are of primary interest in this paper.
Each category is defined in terms of a semidualizing module, 
the study of which
was initiated independently (with different names)
by Foxby~\cite{foxby:gmarm}, Golod~\cite{golod:gdagpi},
and Vasconcelos~\cite{vasconcelos:dtmc}.

\begin{notationdefinition} \label{notation08a}
An $R$-module $C$ is \emph{semidualizing} if it satisfies the following conditions:
\begin{enumerate}[\quad(1)]
\item $C$ admits a (possibly unbounded) resolution by finite rank free $R$-modules,
\item the natural homothety map $R\to\Hom_R(C,C)$ is an isomorphism, and
\item $\ext_R^{\geq 1}(C,C)=0$.
\end{enumerate}
A finitely generated projective  $R$-module of rank 1 is semidualizing.
If $R$ is Cohen-Macaulay, then $D$  is \emph{dualizing}
if it is semidualizing and $\id_R(D)$ is finite.
If $C$ is semidualizing and $D$ is dualizing,
then~\cite[(2.12)]{christensen:scatac} says that
the $R$-module $\da{C}=\hom_R(C,D)$ is semidualizing,
$\aext^{\geq 1}(C,D)=0$ and $C^{\dagger\dagger}\cong C$;
see also~\cite[(4.11)]{vasconcelos:dtmc}.
\end{notationdefinition}

Based on the work of Enochs and Jenda~\cite{enochs:gipm},
the following notions were introduced and studied in this generality by
Holm and J\o rgensen~\cite{holm:smarghd}
and White~\cite{white:gpdrsm}.

\begin{defn} \label{notation08b}
Let $C$ be a semidualizing $R$-module.
An $R$-module is \emph{$C$-projective}
(resp., \emph{$C$-flat}
or \emph{$C$-injective})
if it is isomorphic to a module of the form $P\otimes_R C$ for some projective 
$R$-module $P$
(resp.,  $F\otimes_R C$ for some flat $R$-module $F$
or $\Hom_R(C,I)$ for some injective $R$-module $I$).
We let $\catpc$, $\catfc$ and $\catic$ denote the categories of 
$C$-projective, $C$-flat
and $C$-injective $R$-modules, respectively.

A \emph{complete $\catp\catpc$-resolution} is a complex $X$ of $R$-modules 
satisfying the following:
\begin{enumerate}[\quad(1)]
\item $X$ is exact and $\Hom_R(-,\catpc)$-exact, and
\item $X_i$ is projective when  $i\geq 0$ and $X_i$ is  $C$-projective when $i< 0$.
\end{enumerate}
An $R$-module $G$ is \emph{$\text{G}_C$-projective} if there
exists a complete $\catp\catpc$-resolution $X$ such that $G\cong\coker(\partial^X_1)$,
in which case $X$ is a \emph{complete $\catp\catpc$-resolution of $G$}.  We let
$\catgpc$ denote the subcategory of $\text{G}_C$-projective $R$-modules
and set $\catg\catp=\catg\catp_R$.
Projective $R$-modules and $C$-projective $R$-modules are $\text{G}_C$-projective.

The terms \emph{complete $\catic\cati$-coresolution}
and \emph{$\text{G}_C$-injective} are defined dually, and
$\catgic$ is the subcategory of $\text{G}_C$-injective $R$-modules.
An $R$-module that is injective or $C$-injective is $\text{G}_C$-injective.

Assume that $R$ is noetherian.
A \emph{complete $\catf\catfc$-resolution} is a complex $X$ of $R$-modules 
satisfying the following conditions:
\begin{enumerate}[\quad(1)]
\item $X$ is exact and $-\otimes_R \catic$-exact, and
\item $X_i$ is flat when  $i\geq 0$ and $X_i$ is  $C$-flat when $i< 0$.
\end{enumerate}
An $R$-module $G$ is \emph{$\text{G}_C$-flat} if there
exists a complete $\catf\catfc$-resolution $X$ such that $G\cong\coker(\partial^X_1)$,
in which case $X$ is a \emph{complete $\catf\catfc$-resolution of $G$}.  We let
$\catgfc$ denote the subcategory of $\text{G}_C$-flat $R$-modules
and set $\catg\catf=\catg\catf_R$.
Flat $R$-modules
(hence, projective $R$-modules) and $C$-flat $R$-modules are $\text{G}_C$-flat.
\end{defn}

The $\text{G}_C$-flats are only used in this paper as a tool for 
verifying certain relations between
$\text{G}_C$-projectives and $\text{G}_C$-injectives.
These relations are contained in the next
result which is essentially an assemblage of facts from~\cite{holm:smarghd}.

\begin{lem} \label{gflat}
Assume that $R$ is noetherian.
Let $C$, $E$ and $M$ be $R$-modules with
$C$ semidualizing and $E$ faithfully injective.
\begin{enumerate}[\quad\rm(a)]
\item \label{gflatitem1}
There is an inequality $\gfcpd_R(M)\leq\gpcpd_R(M)$,
and so $\catgpc\subseteq\catgfc$. 
\item \label{gflatitem2}
If $M$ is $\text{G}_C$-flat, then $\ahom(M,E)$ is $\text{G}_C$-injective.
\item \label{gflatitem3}
If $R$ has finite Krull dimension, then the quantities
$\gicid_R(\ahom(M,E))$, $\gpcpd_R(M)$ and $\gfcpd_R(M)$
are simultaneously finite.
\end{enumerate}
\end{lem}

\begin{proof}
\eqref{gflatitem1}
Let $R\ltimes C$ denote the trivial extension of $R$ by $C$ and 
view $M$ as an $R\ltimes C$-module via the
natural surjection $R\ltimes C\to R$.  
In the next sequence 
\begin{align*}
\gfcpd_R(M)
=\gfd_{R\ltimes C}(M)
\leq\gpd_{R\ltimes C}(M)
=\gpcpd_R(M)
\end{align*}
the equalities are from~\cite[(2.16)]{holm:smarghd}
and the inequality is from~\cite[(5.1.4)]{christensen:gd}.

\eqref{gflatitem2}
Assume $M\in\catgfc$.
From~\cite[(2.16)]{holm:smarghd} we know that $G$ is 
Gorenstein flat over $R\ltimes C$, and so~\cite[(2.15)]{holm:smarghd}
implies that $\ahom(M,E)$ is Gorenstein injective over $R\ltimes C$.
An application of~\cite[(2.13.1)]{holm:smarghd} implies
$\ahom(M,E)\in \catgic$.

\eqref{gflatitem3}
Set $\mat{(-)}=\ahom(-,E)$.
Using~\cite[(2.1),(2.16)]{holm:smarghd}
we have equalities
\begin{gather*}
\gicid_R(\mat{M})
=\gid_{R\ltimes C}(\mat{M})
=\gfd_{R\ltimes C}(M)
=\gfcpd_R(M) \\
\gpbpd_R(C)=\gpd_{R\ltimes B}(C)
\end{gather*}
From~\cite[(3.4)]{esmkhani:ghdac} we conclude that
$\gfd_{R\ltimes C}(M)$ and $\gpd_{R\ltimes C}(M)$
are simultaneously finite, and hence so are the six displayed
quantities.
\end{proof}

The following equalities are taken from~\cite[(2.11)]{takahashi:hasm}.

\begin{fact}\label{dimension}
Let $C$ and $M$ be $R$-modules with
$C$ semidualizing.
\begin{enumerate}[\quad\rm(a)]
\item\label{dimensionitem1}
$\pd_R(M)  =\pcdim_R(C\otimes_R M)$ and $\pcdim_R(M)  = \pd_R(\Hom_R(C,M))$.
\item\label{dimensionitem2}
$\icdim_R(M)  =\id_R(C\otimes_R M)$ and $\id_R(M)  = \icdim_R(\Hom_R(C,M))$.
\end{enumerate}
\end{fact}

\begin{defn} \label{appx}
Let $M$ and $N$ be $R$-modules such that $\gpcpd_R(M)<\infty$ 
and $\gicid_R(N)<\infty$.
From~\cite[(4.6) and its dual]{white:gpdrsm} 
there are exact sequences
\begin{align*}
0\to & K \to G_0 \to M \to 0
& 0\to & N\to H_0\to L \to 0
\end{align*}
such that 
$\pcpd_R(K)$ and $\icid_R(L)$ are finite, 
$G_0$ is $\text{G}_C$-projective,
and $H_0$ is $\text{G}_C$-injective.
The first exact sequence is called a \emph{$\catgpc$-approximation} of $M$,
and the second one is called a \emph{$\catgic$-coapproximation} of $N$.

Augmenting the $\catgpc$-approximation with a bounded $\catpc$-resolution
of $K$ yields a bounded $\catgpc$-resolution
$G\xra{\simeq}M$ such that $G_n\in\catpc$ for each $n\geq 1$.
Such a resolution is called a \emph{bounded strict $\catgpc$-resolution}.
Dually, $N$ admits a \emph{bounded strict $\catgic$-coresolution}.
\end{defn}

The next definition was first introduced by Auslander and 
Bridger~\cite{auslander:adgeteac,auslander:smt}
in the case $C=R$, and in this generality
by Golod~\cite{golod:gdagpi} and Vasconcelos~\cite{vasconcelos:dtmc}.

\begin{defn} \label{notation08c}
Assume that $R$ is noetherian,  and let $C$ be a
semidualizing $R$-module.  A finitely generated $R$-module
$H$ is \emph{totally $C$-reflexive} if 
\begin{enumerate}[\quad(1)]
\item $\ext_R^{\geq 1}(H,C)=0=\ext_R^{\geq 1}(\Hom_R(H,C),C)$, and 
\item the natural biduality map $H\to\Hom_R(\Hom_R(H,C),C)$ is an isomorphism.
\end{enumerate}
A finitely generated
module that is projective or $C$-projective is totally $C$-reflexive.
Let
$\catgc$ denote the subcategory of totally $C$-reflexive $R$-modules
and set $\catg=\catg_R$.
\end{defn}

\begin{fact} \label{useful}
The category $\catpc$ is an injective cogenerator for 
$\catgpc$ by~\cite[(2.5),(2.13)]{holm:smarghd} 
and~\cite[(3.2),(3.9)]{white:gpdrsm},
and $\catic$ is a projective generator for 
$\catgic$ by~\cite[(2.6),(2.13)]{holm:smarghd} 
and results dual to~\cite[(3.2),(3.9)]{white:gpdrsm}.
Lemma~\ref{gencat01} yields the relations
$\catgpc\perp\finrescatpc$ and $\fincorescatic\perp \catgic$.
From~\cite[(5.6)]{holm:smarghd} there is an equality
$\propcorescatgic=\catm$.

Let $M$ and $N$ be $R$-modules such that $\gpcpd_R(M)<\infty$ 
and $\gicid_R(N)<\infty$.
The proof of~\cite[(4.6)]{white:gpdrsm} shows that $M$ admits a
bounded strict $\catgpc$-resolution such that
$G_n=0$ for each $n>\gpcpd_R(M)$, and~\cite[(4.4)]{white:gpdrsm} shows 
that every bounded strict $\catgpc$-resolution of $M$ is 
is $\catgpc$-proper and hence $\catpc$-proper. In particular, 
every bounded $\catpc$-resolution is $\catgpc$-proper
and every $\catgpc$-approximation is $\ahom(\catgpc,-)$-exact.
Dually, $N$ admits a bounded strict $\catgic$-coresolution
$N\xra{\simeq}H$ such that $H_{-n}=0$ for each $n>\gicid_R(M)$,
every bounded strict $\catgic$-coresolution
is $\catgic$-proper, and
every bounded $\catic$-coresolution  is $\catgic$-proper.

Assuming that $R$ is noetherian, the equality $\catgc=\catgpc^f$ is
by~\cite[(5.4)]{white:gpdrsm}, and $\catpc^f$ is an injective cogenerator for 
$\catgc$ by~\cite[(3.9),(5.3),(5.4)]{white:gpdrsm}.  
\end{fact}

\begin{notation} \label{notn1}
We simplify our notation for the relative cohomologies 
\begin{align*}
\pcext^n(-,-)&=\ext^n_{\catpc R}(-,-)
&\gpcext^n(-,-)&=\ext^n_{\catgpc R}(-,-) \\
\icext^n(-,-)&=\ext^n_{R\,\catic}(-,-) 
&\gicext^n(-,-)&=\ext^n_{R\,\catgic}(-,-)
\intertext{and for the various connecting maps from Definition~\ref{rel02a}}
\gpcpccomp^n&=\vartheta_{\catgpc\catpc R}^n
& \giciccomp^n&=\vartheta_{R\catgic\catic}^n \\
\gpccomp^n&=\varkappa_{\catgpc R}^n
& \giccomp^n&= \varkappa_{R\catgic}^n \\
\pccomp^n&=\varkappa_{\catpc R}^n
& \iccomp^n&= \varkappa_{R\catic}^n. 
\end{align*}
Fact~\ref{useful} implies that each bifunctor
$\gicext^n(-,-)$ is defined on $\catm\times\catm$.
\end{notation}

The next properties are from~\cite[(4.1)]{takahashi:hasm}.

\begin{fact}\label{relabs0}
Let $C$, $M$ and $N$ be $R$-modules with
$C$ semidualizing.
\begin{enumerate}[\quad\rm(a)]
\item\label{relabs0item1} 
If $M\in\proprescatpc$, then there is an isomorphism for each $n$
$$\ext^n_{\PP_C}(M,N)\cong\ext^n_R(\Hom_R(C,M),\Hom_R(C,N)).$$
\item\label{relabs0item2} 
If $N\in\propcorescatic$, then there is an isomorphism for each $n$
$$\ext^n_{\I_C}(M,N)\cong\ext^n_R(C\otimes_R M,C\otimes_R N).$$
\end{enumerate}
\end{fact}

The following is for use in Propositions~\ref{notiso01} and~\ref{notiso02}.

\begin{lem} \label{useful2}
Assume that $R$ is noetherian and let $C$, $M$ and $N$ be finitely generated
$R$-modules with $C$  semidualizing.
\begin{enumerate}[\quad\rm(a)]
\item \label{useful2item0}
A $\catpc$-resolution
$X\xra{\simeq}M$ is proper if and only if  $\ahom(C,X^+)$ is 
exact.
\item \label{useful2item1}
Assume $M\in\proprescatpc$.
Then $\pcext^n(M,N)\cong\ext^n_{\smash[t]{\catpc^f}}(M,N)$ is finitely generated and
$\supp(\pcext^n(M,N))\subseteq\supp(M)\cap\supp(N)$
for each $n$.  Hence, if $\supp(M)\cap\supp(N)\subseteq\mspec(R)$,
then $\pcext^n(M,N)$ has finite length.
\item \label{useful2item2}
Assume $\gpcpd_R(M)<\infty$.
Then $\gpcext^n\!(M,N)\cong\ext^n_{\catgc}(M,N)$ is finitely generated and
$\supp(\gpcext^n(M,N))\subseteq\supp(M)\cap\supp(N)$
for each $n$. If $\supp(M)\cap\supp(N)\subseteq\mspec(R)$,
then $\gpcext^n(M,N)$ has finite length.
\end{enumerate}
\end{lem}

\begin{proof}
\eqref{useful2item0}
This is immediate from Hom-tensor adjointness.

\eqref{useful2item1}
As $M$ admits
a proper $\catpc$-resolution, we know from~\cite[(2.4.c)]{white:gpdrsm}
that $M$ admits a  $\catpc^f$-resolution $X\xra{\simeq}M$ that is
$\catpc$-proper.
This yields the isomorphism $\pcext^n(M,N)\cong\ext^n_{\smash[t]{\catpc^f}}(M,N)$
and the finite generation of $\ext^n_{\smash[t]{\catpc^f}}(M,N)$.

Fix a prime $\p\in \spec(R)$.  
It is straightforward to check that
the localization $C_{\p}$ is $R_{\p}$-semidualizing and,
using part~\eqref{useful2item0}, that the $\catp_{C_{\p}}$-resolution
$X_{\p}\xra{\simeq}M_{\p}$ is proper.
This yields an isomorphism
$\pcext^n(M,N)_{\p} \cong\ext^n_{\catp_{C_{\p}}}(M_{\p},N_{\p})$
for each integer $n$.  If $\p\not\in\supp(M)$, then the complex
$X_{\p}$ is exact and hence split-exact by~\cite[(2.5)]{white:gpdrsm};
it follows easily that 
$\pcext^n(M,N)_{\p} \cong\ext^n_{\catp_{C_{\p}}}(M_{\p},N_{\p})=0$.
If $\p\not\in\supp(N)$, then one derives the same vanishing.

\eqref{useful2item2}
Using~\cite[(5.6)]{white:gpdrsm}, the assumption
$\gpcpd_R(M)<\infty$ implies that $M$ admits a bounded strict
$\catgpc$-resolution $G$ such that $G_n$ is finitely generated for each $n\geq 0$.  
It follows that the localized complex
$G_{\p}$ is a bounded strict $\catg\catp_{C_{\p}}$-resolution of $M_{\p}$
for each $\p\in\spec(R)$, and 
hence it is a proper $\catg\catp_{C_{\p}}$-resolution 
by Fact~\ref{useful}.
Furthermore, if $M_{\p}=0$, then $G_{\p}$ is a bounded augmented
$\catp_{C_{\p}}$-resolution of $(G_0)_{\p}$, and it follows from~\cite[(2.5)]{white:gpdrsm}
that $G_{\p}$ is split-exact.
The proof now concludes as in part~\eqref{useful2item1}.
\end{proof}

Over a noetherian ring, the next categories were introduced by 
Avramov and Foxby~\cite{avramov:rhafgd}
when $C$ is dualizing, and 
by Christensen~\cite{christensen:scatac} for arbitrary $C$.
(Note that these works (and others) use the notation $\catac$
and $\catbc$ for certain categories of complexes, while our 
categories consist precisely of the modules in these other categories.)
In the non-noetherian setting, these definitions are from~\cite{holm:fear,white:gpdrsm}.

\begin{defn} \label{notation08d}
Let $C$ be a
semidualizing $R$-module.  
The \emph{Auslander class} of $C$ is the subcategory $\catac$
of $R$-modules $M$ such that 
\begin{enumerate}[\quad(1)]
\item $\tor^R_{\geq 1}(C,M)=0=\ext_R^{\geq 1}(C,C\otimes_R M)$, and
\item The natural map $M\to\Hom_R(C,C\otimes_R M)$ is an isomorphism.
\end{enumerate}
The \emph{Bass class} of $C$ is the subcategory $\catbc$
of $R$-modules $M$ such that 
\begin{enumerate}[\quad(1)]
\item $\ext_R^{\geq 1}(C,M)=0=\tor^R_{\geq 1}(C,\Hom_R(C,M))$, and 
\item The natural evaluation map $C\otimes_R\Hom_R(C,M)\to M$ is an isomorphism.
\end{enumerate}
\end{defn}

\begin{fact} \label{GP}
Let $C$ be a semidualizing $R$-module, and
set $\opg(\catpc)=\catgpc\cap\catbc$ and 
$\opg(\catic)=\catgic\cap\catac$.
The category $\catpc$
is an injective cogenerator and
a projective generator for $\opg(\catpc)$ by~\cite[(5.3)]{sather:sgc}.
Dually, the category $\catic$
is an injective cogenerator and
a projective generator for $\opg(\catic)$ by~\cite[(5.4)]{sather:sgc}.

Assume that $R$ is noetherian and set 
$\opg(\catpc^f)=\catgc\cap\catbc=\opg(\catpc)^f$.
The category $\catpc^f$
is an injective cogenerator and
a projective generator for $\opg(\catpc^f)$ by~\cite[(5.5)]{sather:sgc}.
If $R$ is Cohen-Macaulay with a dualizing module, 
then
there are containments 
$\catgpc\subseteq\catacd$ and $\catgic\subseteq\catbcd$
by~\cite[(4.6)]{holm:smarghd}, and 
we conclude
$\opg(\catpc^f)\subseteq\opg(\catpc)\subseteq\cata_{\da{C}}\cap\catbc$
and
$\opg(\catic)\subseteq\catb_{\da{C}}\cap\catac$.
\end{fact}

\begin{fact}\label{projac}
Let $C$ be a
semidualizing $R$-module.  
If any two $R$-modules in a short exact
sequence are in $\cata_C$, respectively $\catb_C$, then so is the
third; see~\cite[(6.7)]{holm:fear}.
The class $\cata_C$ contains all modules of finite
projective dimension and those of finite 
$\catic$-injective dimension, and the class $\catb_C$ contains all modules of
finite injective dimension and those of finite 
$\catpc$-projective dimension
by~\cite[(6.4),(6.6)]{holm:fear}.
If $M$ is in $\catbc$, then $M$ admits a 
proper $\catpc$-resolution; if $M$ is in $\catac$, then
$M$ admits a 
proper $\catic$-injective coresolution; see~\cite[(2.4)]{takahashi:hasm}.

Using the containment $\catggpc\subseteq\catbc$ and 
a $\catgpc$-approximation,
one checks readily that $\ggpcpd_R(M)$ is finite if and only if
$\gpcpd_R(M)$ is finite and $M$ is in $\catbc$.
Consequently, if $\ggpcpd_R(M)$ is finite
(e.g., if $\pcpd_R(M)$ is finite), then 
$M$ admits a 
proper $\catpc$-resolution, and 
$\ggpcpd_R(M)=\gpcpd_R(M)$.

Dually, 
$\ggicid_R(M)$ is finite if and only if
$\gicid_R(M)$ is finite and $M$ is in $\catac$.
If $\ggicid_R(M)$ is finite
(e.g., if $\icid_R(M)$ is finite), then 
$M$ admits a 
proper $\catic$-coresolution, and 
$\ggicid_R(M)=\gicid_R(M)$.

If $R$ is Cohen-Macaulay with a dualizing module, 
then Fact~\ref{GP} yields
$$\finrescatpc\subseteq\finrescatggpc\subseteq\catacd\cap\catbc
\supseteq \fincorescatggicd \supseteq \fincorescaticd.$$
\end{fact}

The following relations between semidualizing modules
are for use in Section~\ref{sec08}.

\begin{lem} \label{cincatbc}
Assume that $R$ is noetherian, and let $B$ and $C$ be semidualizing $R$-modules.
The following conditions are equivalent.
\begin{enumerate}[\quad\rm(i)]
\item \label{cincatbcitem1}
$\gpcpd_R(B)$ is finite.
\item \label{cincatbcitem2}
$B$ is totally $C$-reflexive.
\item \label{cincatbcitem4}
$\ext^{\geq 1}_R(B,C)=0$ and $\ahom(B,C)$ is $R$-semidualizing.
\item \label{cincatbcitem3}
$C$ is in $\catbb$.
\end{enumerate}
\end{lem}

\begin{proof}
\eqref{cincatbcitem1}$\implies$\eqref{cincatbcitem2}
If $\gpcpd_R(B)$ is finite then $\gcdim{C}_R(B)<\infty$
and so~\cite[(3.1)]{frankild:rrhffd} provides the
equality $\gcdim{C}_R(B)=0$ and hence the desired conclusion.

\eqref{cincatbcitem2}$\implies$\eqref{cincatbcitem4}
This is in~\cite[(2.11)]{christensen:scatac}.

\eqref{cincatbcitem4}$\implies$\eqref{cincatbcitem3}
Let $P\xra{\simeq} B$ and $C\xra{\simeq} I$ be projective and injective resolutions, respectively.
The condition $\aext^{\geq 1}(B,C)=0$ implies that
$\ahom(P,I)$ is an injective resolution of $\ahom(B,C)$.
Consider the next commutative diagram
\begin{equation} \label{bassdiag}
\begin{split}
\xymatrix{
R \ar[rr]^-{X'} \ar[d]_{X}^{\simeq} && \ahom(\ahom(P,I),\ahom(P,I)) \ar[d]_{\Phi}^{\cong} \\
\ahom(I,I) \ar[rr]^-{\ahom(\Omega,I)} && \ahom(P\otimes_R\ahom(P,I),I)
}
\end{split}
\end{equation}
where $X$ and $X'$ are the homothety homomorphisms,
$\Phi$ is Hom-tensor adjunction,
and $\Omega$ is tensor-evaluation.
Our assumptions imply that $X'$ is a quasiisomorphism,
and so the same is true of $\ahom(\Omega,I)$.  Using~\cite[(A.8.11)]{christensen:gd}
we conclude that $\Omega$ is also a quasiisomorphism; this uses the equality
$\supp_R(C)=\spec(R)$ which holds because $C$ is semidualizing.
In particular, we have
$$\tor^R_{n}(B,\ahom(B,C))
\cong\HH_n(\ahom(P,I)\otimes_R P)
\cong\HH_n(C)
$$
which is 0 when $n\geq 1$.  
The isomorphism $\HH_0(\Omega)$ is exactly the natural evaluation map
$B\otimes_R\ahom(B,C)\to C$, and so we have $C\in\catbb$.

\eqref{cincatbcitem3}$\implies$\eqref{cincatbcitem4}
Assume that $C$ is in $\catbb$ and employ the notation from the
previous paragraph.  
It follows that the morphism $\Omega$ is a quasiisomorphism,
and hence so is $X'$.  This implies that $\ahom(B,C)$ is semidualizing.
The Bass class conditions then conspire with~\cite[(3.1.c)]{frankild:rrhffd}
to imply $\gpcpd_R(B)<\infty$.
\end{proof}

\begin{fact} \label{compare01}
If $B$ and $C$ be semidualizing $R$-modules such that $\gpcpd_R(B)$ is finite, then
there is a containment $\catpb\subseteq\catgpc$, and
$C$ admits a proper $\catpb$-resolution by
Fact~\ref{projac} and Lemma~\ref{cincatbc}.
For example, 
the semidualizing module $B=R$ is always totally $C$-reflexive; if $R$ is Cohen-Macaulay
and $C$ is dualizing, then $B$ is totally $C$-reflexive. 
For discussions of methods for generating 
other nonisomorphic semidualizing modules
$B$ and $C$ such that $\gpcpd_R(B)<\infty$, the interested reader is encouraged
to peruse~\cite{frankild:rrhffd,frankild:sdcms,sather:divisor}.
\end{fact}

\begin{lem} \label{cincatbc'}
Assume that $(R,\m,k)$ is local and let $B$ and $C$ be semidualizing $R$-modules
with $\gpcpd_R(B)<\infty$.  Let $E$ be the $R$-injective hull of $k$,
and set $\mat{(-)}=\ahom(-,E)$.
The following conditions are equivalent.
\begin{align*}
\mathrm{(i)} & \,\, B\cong C.
& \mathrm{(v)} & \,\, \pcpd_R(B)<\infty. \\
 \mathrm{(ii)} & \,\,\pbpd_R(C)<\infty.
& \mathrm{(vi)} & \,\, \ibid_R(\mat{C})<\infty.  \\
\mathrm{(iii)} & \,\, \gpbpd_R(C)<\infty.
&  \mathrm{(vii)} & \,\, \gibid_R(\mat{C})<\infty.  \\
\mathrm{(iv)} & \,\, \pd_R(\ahom(B,C))<\infty. 
& \mathrm{(viii)} & \,\, \icid_R(\mat{B})<\infty. 
\end{align*}
\end{lem}

\begin{proof}
The implication
(i)$\implies$(n) is straightforward for 
$n=\text{i},\ldots,\text{viii}$,
as are
(ii)$\implies$(iii)
and 
(vi)$\implies$(vii).
The implication
(iii)$\implies$(i)
is in~\cite[(5.3)]{takahashi:hiatsb},
and (ii)$\iff$(iv)
is from Fact~\ref{dimension}\eqref{dimensionitem1}, while
(vii)$\iff$(iii)
is in Lemma~\ref{gflat}\eqref{gflatitem3}.

(v)$\implies$(iii)
If $\pcpd_R(B)<\infty$, then
$B$ is in $\catbc$ and so
Lemma~\ref{cincatbc} implies 
$\gpbpd_R(C)<\infty$.

(viii)$\implies$(v)
Assume $\icid_R(\mat{B})<\infty$. 
Hom-evaluation yields an isomorphism
$$\mat{\ahom(C,B)}\cong C\otimes_R(\mat{B})$$
and hence the first equality in the following sequence
$$\id_R(\mat{\ahom(C,B)})
=\id_R(C\otimes_R(\mat{B}))
=\icid_R(\mat{B})<\infty.$$
The second equality is by Fact~\ref{dimension}\eqref{dimensionitem2}.
It follows that $\ahom(C,B)$ has finite projective dimension
and so Fact~\ref{dimension}\eqref{dimensionitem1} implies
$\pcpd_R(B)<\infty$.
\end{proof}

\section{Comparison Isomorphisms} \label{sec5}

The results of this section document 
situations where different relative cohomology theories agree.
The notation for the comparison homomorphisms is given in~\ref{notn1}.
The next result is a more precise version of~\cite[(4.2)]{takahashi:hasm}.

\begin{prop}\label{relabs}
Let $C$, $M$ and $N$ be  $R$-modules
with $C$ semidualizing.
\begin{enumerate}[\quad\rm(a)]
\item \label{bcabsitem1}
If $M$ and $N$ are in $\catb_C$, then 
the natural map 
$$\pccomp^n(M,N)\colon\ext^n_{\PP_C}(M,N)\to\ext^n_R(M,N)$$
is an isomorphism for each integer $n$.
\item\label{bcabsitem2}
If $M$ and $N$ are in $\cata_C$, then 
the natural map 
$$\iccomp^n(M,N)\colon\icext^n(M,N)\to\ext^n_R(M,N)$$
is an isomorphism for each integer $n$.
\end{enumerate}
\end{prop}

\begin{proof}
\eqref{bcabsitem1}
Let $P\xra{\gamma}\hom_R(C,M)$ and $P'\xra{\gamma'}C$ be projective
resolutions.
Because $M$ is in $\catbc$, we have $\tor^R_{\geq 1}(C,\hom_R(C,M))=0$
and so the complex $P'\otimes_R P$ is a projective resolution of
$C\otimes_R\hom_R(C,M)\cong M$, 
and the complex $C\otimes_R P$ is a $\catpc$-resolution of
$M$.  The  following diagram commutes
$$\xymatrix{
P'\otimes_R P \ar[r]^-{\gamma'\otimes_R P} \ar[d]_{\simeq} & C\otimes_R P \ar[d]^{\simeq} \\
M \ar[r]^-{\id_M} & M
}$$
and so it suffices to show that the induced map
$$\ahom(P'\otimes_R P,N)\xra{\ahom(\gamma'\otimes_R P,N)}
\ahom(C\otimes_R P,N)$$
is a quasiisomorphism.
The following standard isomorphisms
\begin{align*}
\cone(\ahom(\gamma'\otimes_R P,N))
&\cong \shift\ahom(\cone(\gamma'\otimes_R P),N)\\
&\cong \shift\ahom(\cone(\gamma')\otimes_R P,N)
\end{align*}
imply that it suffices to show that the complex 
$\ahom(\cone(\gamma')\otimes_R P,N)$ is exact.  

Observe that $\cone(\gamma')$ is exact and bounded below and
each module $\cone(\gamma')_n$ is a direct sum of a projective
$R$-module and a $C$-projective $R$-module.  Since $N$ is in $\catbc$,
we know that $\ext^{\geq 1}_R(C,N)=0$, and it follows that
$\ext^{\geq 1}_R(Q\otimes_RC,N)=0$ for each projective $R$-module $Q$.
Since we also have $\ext^{\geq 1}_R(Q,N)=0$, it follows that
$\ext^{\geq 1}_R(\cone(\gamma')_n,N)=0$ for each $n$.
Breaking up  $\cone(\gamma')$ into short exact sequences
and applying $\ahom(-,N)$ to each piece yields the desired conclusion.

The proof of~\eqref{bcabsitem2} is dual.
\end{proof}

The next result compares to~\cite[(4.2.3)]{avramov:aratc}.

\begin{prop} \label{detect01}
Let $C$ and $M$ be  $R$-modules with $C$ semidualizing.
\begin{enumerate}[\quad\rm(a)]
\item \label{detect01item1}
If $\pcpd_R(M)$ is finite, then 
the following natural transormations are isomorphisms for each $n$
$$\gpcpccomp^n(M,-)\colon\gpcext^n(M,-)\xra{\cong}\pcext^n(M,-)$$ 
and so
$\gpcext^n(M,-)=0$ for each $n>\pcpd(M)$.
\item \label{detect01item2}
If $\icid_R(M)$ is finite, then 
the following natural transormations are isomorphisms for each $n$
$$\giciccomp^n(-,M)\colon\gicext^n(-,M)\xra{\cong}\icext^n(-,M)$$ 
and so $\gicext^n(-,M)=0$ for each $n>\icid(M)$.
\end{enumerate}
\end{prop}

\begin{proof}
We prove part~\eqref{detect01item1};  the proof of~\eqref{detect01item2}
is dual.
Let $W\xra{\simeq}M$ be a $\catpc$-resolution such that $W_n=0$ for each $n>\pcpd_R(M)$.
The resolution $W$ is  $\catgpc$-proper and $\catpc$-proper 
by Fact~\ref{useful},
so both $\gpcext^n(M,-)$ and $\pcext^n(M,-)$ are defined.  Further,
in the notation of Definition~\ref{rel02a}, we can take
$\ol{\id_M}=\id_W$, and so  the natural isomorphisms follow from the next equalities
$$\xwacomp^n(M,-)
=\HH_{-n}(\ahom(\id_W,-))
=\id_{\HH_{-n}(\ahom(W,-))}.$$
The vanishing conclusion follows readily since
$W_n=0$ for each $n>\wpd(M)$.
\end{proof}

The next lemma is a tool for the proofs of Propositions~\ref{detect01a} 
and~\ref{detect01b}.  Note that we do not assume that the complexes satisfy
any properness conditions.

\begin{lem} \label{quisos}
Let $C$, $M$, and $N$ be $R$-modules with $C$ semidualizing.
\begin{enumerate}[\quad\rm(a)]
\item \label{quisositem1}
Let  $\alpha\colon G\to G'$ be a quasiisomorphism between
bounded below complexes in $\catgpc$.
If $\pcpd(N)<\infty$,
then the morphism $\ahom(\alpha,N)\colon\ahom(G',N)\to\ahom(G,N)$
is a quasiisomorphism.
\item \label{quisositem2}
Let  $\beta\colon H\to H'$ be a quasiisomorphism between
bounded above complexes in $\catgic$.
If $\icid(M)<\infty$,
then the morphism $\ahom(M,\beta)\colon\ahom(M,H)\to\ahom(M,H')$
is a quasiisomorphism.
\end{enumerate}
\end{lem}

\begin{proof}
We prove part~\eqref{quisositem1}; the proof of part~\eqref{quisositem2} is dual.

It suffices to show that  $\cone(\ahom(\alpha,N))$ is exact.
From the next isomorphism 
$$\cone(\ahom(\alpha,N))\cong\shift\ahom(\cone(\alpha),N)$$
we need to show that $\ahom(\cone(\alpha),N)$ is exact.  Note that
$\cone(\alpha)$ is an exact, bounded below complex in $\catgpc$.  
Set  $M_j=\ker(\partial^{\cone(\alpha)}_j)$ for each integer $j$, 
and note $M_{j-1}\in\catgpc$ for $j \ll 0$.
Consider the exact sequences
\begin{equation} \label{exact01'} \tag{$\ast_j$}
0\to M_j\to \cone(\alpha)_j\to M_{j-1}\to 0.
\end{equation}
Lemma~\ref{gencat01} and Fact~\ref{useful}
imply $\catgpc\perp N$.
Hence, induction on $j$ using Lemma~\ref{perp03}\eqref{perp03item1} implies
$\aext^{\geq 1}(M_j,N)=0$ for each $j$ and so each sequence $(\ast_j)$ is
$\ahom(-,N)$-exact.
It follows that $\ahom(\cone(\alpha),N)$ is exact.
\end{proof}

The next two results compare to~\cite[(4.2.4)]{avramov:aratc}.
Notice that the condition $M$ is in $\proprescatgpc\cap \proprescatpc$
of Proposition~\ref{detect01a}\eqref{detect01aitem1}
is satisfied when $\ggpcpd_R(M)<\infty$;
see Fact~\ref{GP}.
Also,  part~\eqref{detect01aitem2} uses the equality 
$\propcorescatgic=\catm$ from Fact~\ref{useful}.

\begin{prop} \label{detect01a}
Let $C$, $M$, and $N$ be $R$-modules with $C$ semidualizing.
\begin{enumerate}[\quad\rm(a)]
\item \label{detect01aitem1}
If $M$ is in $\proprescatgpc\cap \proprescatpc$ and 
$N$ is in $\finrescatpc$, then 
the following natural map is an isomorphism for each $n$
$$\gpcpccomp^n(M,N)\colon\gpcext^n(M,N)\xra{\cong}\pcext^n(M,N).$$ 
\item \label{detect01aitem2} 
If $M$ is in $\fincorescatic$ and 
$N$ is in $\propcorescatic$, then 
the following natural map is an isomorphism for each $n$
$$\giciccomp^n(M,N)\colon\gicext^n(M,N)\xra{\cong}\icext^n(M,N).$$ 
\end{enumerate}
\end{prop}

\begin{proof}
We prove part~\eqref{detect01aitem1}; the proof of part~\eqref{detect01aitem2}
is dual.  

The module
$M$ has a proper $\catpc$-resolution 
$\gamma\colon W\to M$ and a proper $\catgpc$-resolution 
$\gamma'\colon G\to M$.
Lemma~\ref{rel01}\eqref{rel01item1}
yields a quasiisomorphism
$\ol{\id_M}\colon W\to G$ 
such that $\gamma=\gamma'\ol{\id_M}$,
and Lemma~\ref{quisos}\eqref{quisositem1} implies that $\ahom(\ol{\id_M},N)$ is a quasiisomorphism.
The result now follows from the definition of
$\gpcpccomp^n(M,N)$.
\end{proof}

The next result is proved like Proposition~\ref{detect01a} using the containment
$\catp\subseteq\catgpc$.

\begin{prop} \label{detect01b}
Let $C$, $M$, and $N$ be $R$-modules with $C$ semidualizing.
\begin{enumerate}[\quad\rm(a)]
\item \label{detect01bitem1}
If $M$ is in $\proprescatgpc$ and 
$N$ is in $\finrescatpc$, then 
the following natural map is an isomorphism for each $n$
$$\gpccomp^n(M,N)\colon\gpcext^n(M,N)\xra{\cong}\aext^n(M,N).$$ 
\item \label{detect01bitem2}
If $M$ is in $\fincorescatic$, then 
the following natural map is an isomorphism for each $n$
\begin{align*}
&&&&&&&&&&&&&&\giccomp^n(M,N)\colon\gicext^n(M,N)\xra{\cong}\aext^n(M,N).
&&&&&&&&&&&&&& \qed
\end{align*}
\end{enumerate}
\end{prop}

The next four  lemmata are tools for the proofs of Propositions~\ref{detect01a'} 
and~\ref{detect01b'} and for Theorem~\ref{balance02a}.  
Part~\eqref{little2item0} of
the first one is a consequence of
Proposition~\ref{relabs}, and parts~\eqref{little2item1} and~\eqref{little2item2}
follows from part~\eqref{little2item0}.
Note that Fact~\ref{projac}
gives conditions guaranteeing $M,N\in\catb_C\cap\A_{C^{\dagger}}$.

\begin{lem} \label{little1}
Let $R$ be a Cohen-Macaulay ring with dualizing
module, and let $C$, $M$ and $N$ be $R$-modules with
$C$ semidualizing.
\begin{enumerate}[\quad\rm(a)]
\item \label{little2item0}
If $M,N\in\catb_C\cap\A_{C^{\dagger}}$, then the natural map
$$\ext_{\PP_C}^n(M,N)
\xra[\cong]{\pccomp^n(M,N)}
\ext^n_R(M,N)
\xla[\cong]{\icdcomp^n(M,N)}
\ext^n_{\I_{C^\dagger}}(M,N)
$$
is an isomorphism for each $n$.
\item \label{little2item1}
If $M\in\catb_C\cap\A_{C^{\dagger}}$ and $N\in\cati_{\da{C}}$, then 
$\ext_{\PP_C}^{\geq 1}(M,N)=0=
\ext^{\geq 1}_R(M,N)$.
\item \label{little2item2}
If $M\in\catp_C$ and $N\in\catb_C\cap\A_{C^{\dagger}}$, then 
$\ext_{\I_{C^\dagger}}^{\geq 1}(M,N)=0=
\ext^{\geq 1}_R(M,N)$. \qed
\end{enumerate}
\end{lem}

\begin{lem}\label{orthog}
Let $C$ be a semidualizing $R$-module.
One has $\catac\perp\catic$ and $\catpc\perp\catbc$.
If $R$ is Cohen-Macaulay and admits a dualizing
module, then
$\G\PP_C\perp\I_{C^\dagger}$ and $\PP_C\perp\G\I_{C^\dagger}$,
and so $\PP_C\perp\I_{C^\dagger}$.
\end{lem}

\begin{proof}
We verify the first orthogonality condition; the second one is verified similarly,
and the others follow immediately from the containments
$\catpc\subseteq\catgpc\subseteq\catacd$ and
$\caticd\subseteq\catgicd\subseteq\catbc$; see Fact~\ref{projac}.
Let $M\in \catac$ and $N\in \catic\subseteq\catac$.
For each $n\geq 1$, the isomorphism in the following sequence is
in Proposition~\ref{relabs}\eqref{bcabsitem2}
$$\ext^n_R(M,N)\cong\icext^n(M,N)=0$$
and the vanishing holds because $N\in\catic$.
\end{proof}

\begin{lem} \label{moreperp}
If $C$ is a semidualizing $R$-module, then one has
$\catgpc\perp\fincorescati$
and $\finrescatp\perp\catgic$.
\end{lem}

\begin{proof}
We verify the first orthogonality condition; the second one is verified similarly. 
Fix  modules $G_0\in\catgpc$ and $N\in\fincorescati$ and set
$j=\id_R(N)<\infty$. 
For each $n\geq 0$
use the fact  that $\catpc$ is a cogenerator for $\catgpc$
to find exact sequences
\begin{equation} \tag{$\ast_n$}
0\to G_n\to W_n\to G_{n+1}\to 0
\end{equation}
with $G_{n+1}\in\catgpc$ and $W_n\in\catpc$; see Fact~\ref{useful}.
From Fact~\ref{projac} we know $N\in\catbc$ and so Lemma~\ref{orthog}
implies $\catpc\perp N$.  Hence, for $i>0$
 the long exact sequences in
$\aext(-,N)$ associated to $(\ast_n)$ yield the isomorphism in the following sequence
$$\aext^i(G_0,N)\cong\aext^{i+j}(G_j,N)=0$$
while the vanishing holds because $i+j>j=\id_R(N)$.
\end{proof}

\begin{lem} \label{quisos'}
Let $C$, $M$, and $N$ be $R$-modules with $C$ semidualizing.
\begin{enumerate}[\quad\rm(a)]
\item \label{quisos'item1}
Let  $\alpha\colon G\to G'$ be a quasiisomorphism between
bounded below complexes in $\catgpc$.
If $\id(N)<\infty$,
then the morphism $\ahom(\alpha,N)\colon\ahom(G',N)\to\ahom(G,N)$
is a quasiisomorphism.
\item \label{quisos'item2}
Let  $\beta\colon H\to H'$ be a quasiisomorphism between
bounded above complexes in $\catgic$.
If $\pd(M)<\infty$,
then the morphism $\ahom(M,\beta)\colon\ahom(M,H)\to\ahom(M,H')$
is a quasiisomorphism.
\end{enumerate}
\end{lem}

\begin{proof}
We prove part~\eqref{quisos'item1}; the proof of part~\eqref{quisos'item2} is dual.
Set $M_j=\ker(\partial_j^{\cone(\alpha)})$ for each $j$
and consider the following exact sequences
\begin{equation} \label{another}
0\to M_j\to \cone(\alpha)_j\to M_{j-1}\to 0.
\end{equation}
Because $\cone(\alpha)$ is an exact bounded below complex in $\catgpc$,
we know 
$M_j\in\catgpc$ for $j\ll 0$.
From~\cite[(3.8)]{white:gpdrsm} we know that
$\catgpc$ is closed under kernels of epimorphisms,
so an induction argument using~\eqref{another} implies
$M_j\in\catgpc$ for all $j$.
Thus, Lemma~\ref{moreperp}  yields $M_j\perp N$
and $\cone(\alpha)_j\perp N$
for all $j$. 
The long exact sequence in $\ext_R(-,N)$ shows that~\eqref{another}
is $\ahom(-,N)$-exact, and
the conclusion now follows as in the proof of Lemma~\ref{quisos}.
\end{proof}

The next two results follow from Lemma~\ref{quisos'}
in the same way that Propositions~\ref{detect01a}
and~\ref{detect01b} follow from Lemma~\ref{quisos}.

\begin{prop} \label{detect01a'}
Let $C$, $M$, and $N$ be $R$-modules with $C$ semidualizing.
\begin{enumerate}[\quad\rm(a)]
\item \label{detect01a'item1}
If $M$ is in $\proprescatgpc\cap \proprescatpc$ and 
$N$ is in $\fincorescati$, then 
the following natural map is an isomorphism for each $n$
$$\gpcpccomp^n(M,N)\colon\gpcext^n(M,N)\xra{\cong}\pcext^n(M,N).$$ 
\item \label{detect01a'item2} 
If $M$ is in $\finrescatp$ and 
$N$ is in $\propcorescatic$, then 
the following natural map is an isomorphism for each $n$
\begin{align*}
&&&&&&&&&&\hspace{4.5mm}
&&&\giciccomp^n(M,N)\colon\gicext^n(M,N)\xra{\cong}\icext^n(M,N).
&&&&&&&\hspace{1mm}
&&& \qed
\end{align*}
\end{enumerate}
\end{prop}

\begin{prop} \label{detect01b'}
Let $C$, $M$, and $N$ be $R$-modules with $C$ semidualizing.
\begin{enumerate}[\quad\rm(a)]
\item \label{detect01b'item1}
If $M$ is in $\proprescatgpc$ and 
$N$ is in $\fincorescati$, then 
the following natural map is an isomorphism for each $n$
$$\gpccomp^n(M,N)\colon\gpcext^n(M,N)\xra{\cong}\aext^n(M,N).$$ 
\item \label{detect01b'item2}
If $M$ is in $\finrescatp$, then 
the following natural map is an isomorphism for each~$n$
\begin{align*}
&&&&&&&&&\hspace{7mm}
&&&&\giccomp^n(M,N)\colon\gicext^n(M,N)\xra{\cong}\aext^n(M,N).
&&&&&&& 
&&&&& \qed
\end{align*}
\end{enumerate}
\end{prop}

\section{Balance for Relative Cohomology} \label{sec7} 

This section focuses on balance for the functors
$\pcext^n(-,-)$ and $\ibext^n(-,-)$, and 
for $\gpcext^n(-,-)$ and $\gibext^n(-,-)$.

\begin{defn} \label{balance00}
Fix subcategories $\catx'\subseteq \proprescatx$ and $\caty'\subseteq\propcorescaty$.
We say that $\xaext$ and $\ayext$ are \emph{balanced} on
$\catx'\times\caty'$ when the following condition holds:
for each object $M$ in $\catx'$ and $N$ in $\caty'$, if $X\xra{\simeq}M$ is a proper
$\catx$-resolution, and $N\xra{\simeq} Y$ a proper $\caty$-coresolution, then the induced
morphisms of complexes
$$
\ahom(M,Y)\to
\ahom(X,Y)\from
\ahom(X,N)
$$
are quasiisomorphisms.  
\end{defn}

\begin{disc}
When 
$\xaext$ and $\ayext$ are balanced on
$\catx'\times\caty'$, there are isomorphisms
$\xaext^n(M,N)\cong\ayext^n(M,N)$
for all $M\in \catx'$ and $N\in \caty'$ and $n\in\mathbb{Z}$.
\end{disc}

The next example shows that 
the naive version of balance for relative cohomology does not
hold: when $D$ is dualizing, one can have
$\pdext^n(M,N)\not\cong\idext^n(M,N)$ and
$\gpdext^n(M,N)\not\cong\gidext^n(M,N)$, 
even if $\pdpd_R(M)<\infty$ and $\idid_R(N)<\infty$.

\begin{ex}\label{ctrex2}
Let $(R,\m,k)$ be a local non-Gorenstein 
Cohen-Macaulay ring of dimension $d>0$ with dualizing
module $D$.  
Assume  that $R$ is 
Gorenstein on the punctured spectrum.  
From~\cite[(2.12)]{takahashi:hasm} we have 
$\idid_R(R) =\id_R(D)= d$, and
Proposition~\ref{detect01} provides isomorphisms
$$\gpdext^n(D,R)\cong
\ext^n_{\catp_D}(D,R) \qquad \qquad\ext^n_{\cati_D}(D,R)
\cong\gidext^n(D,R)$$ 
for each $n$.
One has $\pdpd_R(D) = 0$ because $D$ is in $\catpd$,
and so $\ext_{\PP_D}^{\geqslant 1}(D,R) = 0$.
Fact~\ref{relabs0}\eqref{relabs0item2}
yields an isomorphism
$$\ext_{\I_D}^n(D,R)\cong \ext_R^n(D\otimes_RD,D)$$ 
for each integer $n$.
To complete the example, we verify
$$
\ext^n_{\catp_D}(D,R)\ncong\ext^n_{\cati_D}(D,R)
$$ 
for some $n\geq 1$.  
From the vanishing $\ext_{\PP_D}^{\geqslant 1}(D,R) = 0$, 
it suffices to find an integer $n\geq 1$ such that
$\ext_R^n(D\otimes_RD,D)\ne 0$.
We utilize the spectral sequence
$$E_2^{p,q}=\ext^p_R(\tor^R_q(D,D),D)\implies\ext^{p+q}_R(D,R).$$
For each  prime ideal $\p\subsetneq\m$, the module $D_{\p}$ is
dualizing for the Gorenstein ring $R_{\p}$,
and so $D_{\p} \cong R_{\p}$; see~\cite[(3.3.7.a)]{bruns:cmr}. 
It follows that $\tor_q^R(D,D)$ has finite length for each
$q>0$, and so $E_2^{p,q}=0$ if $p \ne d$ and $q>0$. For each 
$n\leq d$, this yields 
$$\ext^n_R(D\otimes_R D,D)\cong\ext^n_R(D,R).$$ 
Because $R$ is not Gorenstein, we deduce from~\cite[(2.1)]{avramov:edcrcvct}
that the displayed modules are nonzero for some integer $n$ such that
$1\leq n\leq d$, as desired.
\end{ex}

The following result
contains part of Theorem~\ref{introbalance} from the introduction.

\begin{prop} \label{balance02b}
Assume that $R$ is Cohen-Macaulay ring and admits a dualizing module, and 
let $C$ be a semidualizing $R$-module.
Then $\pcext$ and $\icdext$ are balanced on 
$\finrescatpc\times\fincorescaticd$. In particular, if
$\pcpd_R(M)<\infty$
and $\icdagdim_R(N)<\infty$, then 
there are isomorphisms for each integer $n$
$$\ext_{\PP_C}^n(M,N)\cong\ext^n_{\I_{C^\dagger}}(M,N).$$
\end{prop}

\begin{proof}
Lemma~\ref{little1} implies
$\pcext^{\geq 1}(\finrescatpc,\caticd)=0=
\icext^{\geq 1}(\catpc, \fincorescaticd)$
and so the desired conclusion follows from~\cite[(8.2.14)]{enochs:rha}.
\end{proof}

The next two lemmata are the primary tools for
Theorem~\ref{balance02a}.

\begin{lem} \label{balance01}
Let $R$ be a Cohen-Macaulay ring with dualizing
module. Let $C$, $M$ and $N$ be  $R$-modules
with $C$ semidualizing.
\begin{enumerate}[\quad\rm(a)]
\item \label{balance01item1}
If $N\in\catg\cati_{\da{C}}$ and $\pcdim_R(M)<\infty$,
then $\ext_{\PP_C}^{\geq 1}(M,N)=0$.
\item \label{balance01item2}
If $M\in\catg\catp_{C}$ and $\icdagdim_R(N)<\infty$, 
then $\ext_{\cati_{\da{C}}}^{\geq 1}(M,N)=0$.
\end{enumerate}
\end{lem}

\begin{proof}
We prove part~\eqref{balance01item1};  part~\eqref{balance01item2} is  dual.  
Set $N_0=N$, and for each $n\geq 0$
use the fact  that $\caticd$ is a generator for $\catgicd$  
to find exact sequences 
\begin{equation} \label{another1}
0\to N_{n+1}\to V_n\to N_n\to 0
\end{equation}
with $V_n$ in $\caticd$ and $N_{n+1}$ in $\catgicd$; 
see Fact~\ref{useful}.  Lemma~\ref{orthog} implies
$\catpc\perp\catgicd$, and so
the long-exact sequence in $\ext_R(\catpc,-)$
shows that~\eqref{another1} is  $\ahom(\catpc,-)$-exact.
Fix an integer $j\geq 1$ and set $p=\pcpd(M)$.
Lemma~\ref{little1}\eqref{little2item1} implies $\pcext^{\geq 1}(M,V_n)=0$ for each $n$.
Hence, Lemma~\ref{notation06a}~\eqref{06aitem1} and Remark~\ref{rel03}
yield the isomorphism in the next sequence
$$\pcext^j(M,N)=\pcext^j(M,N_0)\cong\pcext^{j+p}(M,N_{p})=0 $$
where the vanishing is from the (in)equalities $j+p>p=\pcpd_R(M)$.
\end{proof}

\begin{lem} \label{balance04a}
Let $R$ be a  Cohen-Macaulay ring admitting a dualizing
module.  Let $C$, $M$ and $N$ be  $R$-modules
with $C$ semidualizing.
\begin{enumerate}[\quad\rm(a)]
\item \label{balance04aitem1}
Assume that $\gpcpd_R(M)$ is finite and let $G\xra{\alpha} M$
be a proper $\catgpc$-resolution. If $Y$ is
a bounded above complex of objects in $\catgicd$, then the induced map
$\Hom_R(M,Y)\to\Hom_R(G,Y)$ is a quasiisomorphism.
\item \label{balance04aitem2}
Assume that $\gicdid_R(N)$ is finite and let $N\xra{\beta} H$
be a proper $\catgicd$-resolution. If $X$ is
a bounded below complex of objects in $\catgpc$, then the induced map
$\Hom_R(X,N)\to\Hom_R(X,H)$ is a quasiisomorphism.
\end{enumerate}
\end{lem}

\begin{proof}
We proof part~\eqref{balance04aitem1}; the proof
of~\eqref{balance04aitem2} is dual.
To show that the induced map
$\hom_R(\alpha,Y)\colon \Hom_R(M,Y)\to\Hom_R(G,Y)$ is a
quasiisomorphism, it suffices to show that $\cone(\hom_R(\alpha,Y))$
is exact. From the isomorphisms of complexes
$$\cone(\hom_R(\alpha,Y))\cong\shift\hom_R(\cone(\alpha),Y)
\cong\shift\hom_R(G^+,Y)$$ 
and a standard argument, it
suffices to show that $\hom_R(G^+,Y_j)$ is exact for each~$j$.

The module $M$ has a bounded strict $\catgpc$-resolution
$G'\xra{\alpha'} M$ by Fact~\ref{useful}. From Lemma~\ref{rel01}\eqref{rel01item1} 
we conclude  that
$G^+$ and $(G')^+$ are homotopy equivalent, and so the complex
$\hom_R(G^+,Y_j)$ is exact
if and only if $\hom_R((G')^+,Y_j)$ is exact.  Thus, we may replace $G$ with $G'$ 
to assume that $G$ is strict.

For each $n$, set $M_n=\coker(\partial_{n+2}^G)$ and note that $M_{-1}\cong
M$.
For each $n\geq 0$, we have $\pcdim(M_n)<\infty$ and we consider the
following exact sequences
\begin{equation} \label{exseq01a}
0\to M_{n}\xra{\gamma_i} G_{n}\xra{\delta_n} M_{n-1}\to 0.
\end{equation}
It suffices to show that each of these sequences is
$\Hom_R(-,Y_j)$-exact, that is, that the following map is
surjective.
$$\Hom_R(\gamma_n,Y_j)\colon\Hom_R(G_{n},Y_j)\to\Hom_R(M_{n},Y_j)$$
Use the fact  that $\caticd$ is a generator for $\catgicd$  
to find an exact sequence
\begin{equation} \label{exseq02a}
0\to Y'\to V\xra{\tau}Y_j\to 0
\end{equation}
such that $Y'$ is  in $\catgicd$ and $V$ is  in
$\caticd$; see Fact~\ref{useful}. Lemma~\ref{orthog} implies $\catpc\perp\catgicd$ 
and so Lemma~\ref{perp03}\eqref{perp03item2}
guarantees that this sequence
is $\Hom_R(\catpc,-)$-exact.

Fix an element $\lambda\in\Hom_R(M_{n},Y_j)$. The proof will be
complete once we find $f\in\Hom_R(G_n,Y_j)$ such that
$\lambda=f\gamma_n$. The following diagram is our guide
\begin{equation*}
\xymatrix{ && 0\ar[r] & M_n
\ar[r]^{\gamma_n}\ar[d]_>>>>>>{\lambda}\ar@{-->}[ld]_{\sigma}
& G_n \ar[r]\ar@{-->}[lld]_<<<<<<<{\delta}\ar@{-->}[ld]^<<<<{f} & M_{n-1} \ar[r] & 0 \\
0\ar[r] & Y'\ar[r] & V\ar[r]^{\tau} & Y_j \ar[r] & 0 }
\end{equation*}
wherein the top row is~\eqref{exseq01a} and the bottom row
is~\eqref{exseq02a}.

Since the sequence~\eqref{exseq02a} is $\ahom(\catpc,-)$-exact, it gives rise to a
long exact sequence in $\ext_{\PP_C}(M_n,-)$. The vanishing of
$\ext_{\PP_C}^1(M_n,Y')$ 
from Lemma~\ref{balance01}\eqref{balance01item1}
implies that this long exact sequence has
the form
$$0\to \Hom_R(M_n,Y')\to \Hom_R(M_n,V)\xra{\Hom_R(M_n,\tau)} \Hom_R(M_n,Y_j)\to 0.$$
Hence, there exists $\sigma\in\Hom_R(M_n,V)$ such that
$\lambda=\tau\sigma$.

Proposition~\ref{relabs}\eqref{bcabsitem2} implies 
$\ext_R^1(M_{n-1},V)\cong\icdext^1(M,N)=0$, 
so
the long exact sequence in $\ext_R(-,V)$ associated
to~\eqref{exseq01a} has the form
$$0\to \Hom_R(M_{n-1},V)\to \Hom_R(G_n,V)\xra{\Hom_R(\gamma_n,V)} \Hom_R(M_n,V)\to 0.$$
Hence, there exists $\delta\in\Hom_R(G_n,V)$ such that
$\sigma=\delta\gamma_n$.
It follows that
$(\tau\delta)\gamma_n=\tau\sigma=\lambda$
and so $f=\tau\delta\in\Hom_R(G_n,Y_j)$ has the desired property.
\end{proof}

Our main balance result for relative cohomology now follows.
It contains part of Theorem~\ref{introbalance} from the introduction.

\begin{thm} \label{balance02a}
Assume that $R$ is Cohen-Macaulay and admits dualizing module, and 
let $C$ be a semidualizing $R$-module.
Then $\gpcext$ and $\gicdext$ are balanced on 
$\finrescatgpc\times\fincorescatgicd$. In particular, if
$\gpcpd_R(M)<\infty$
and $\gicdagdim_R(N)<\infty$, then 
there are isomorphisms for each integer $n$
$$\ext_{\G\PP_C}^n(M,N)\cong\ext^n_{\G\I_{C^\dagger}}(M,N).$$
\end{thm}

\begin{proof}
From Fact~\ref{useful} we obtain
a bounded proper
$\G\PP_C$-resolution $\alpha\colon G\xra{\simeq}M$ and 
a bounded
proper $\G\I_{C^{\dagger}}$-coresolution $\beta\colon
N\xra{\simeq}Y$. By Lemma~\ref{balance04a}, the
morphisms
$$
\Hom_R(M,Y)\xra{\Hom_R(\alpha,Y)} \Hom_R(G,Y)\xla{\Hom_R(X,\beta)}
\Hom_R(X,N)
$$
are quasiisomorphisms, as desired.
\end{proof}

\section{Distinguishing between Relative Cohomology Theories} \label{sec08}

From Sections~\ref{sec5} and~\ref{sec7} we see that there are numerous situations
where different relative cohomology theories agree.  The purpose of this section is
to show that these theories are almost never \emph{identically} equal.  
Part~\eqref{notiso01item2} of the next result shows 
$\pcext^n(-,-)\ncong\aext^n(-,-)\ncong\gpcext^n(-,-)$.
Part~\eqref{notiso01item3} shows again
$\aext^n(-,-)\ncong\gpcext^n(-,-)$.
Part~\eqref{notiso01item4} shows 
$\pcext^n(-,-)\ncong\gpcext^n(-,-)$.
Lemma~\ref{notisoex} shows how one can construct modules satisfying the
hypotheses of part~\eqref{notiso01item4}.

\begin{prop} \label{notiso01}
Let $(R,\m,k)$ be a local ring and
$C$ a semidualizing $R$-module
such that $C\not\cong R$.
\begin{enumerate}[\quad\rm(a)]
\item \label{notiso01item2}
If $n\geq 1$, then
$\pcext^n(C,k)=0=\gpcext^n(C,k)$ and
$\aext^n(C,k)\neq 0$.
\item \label{notiso01item3}
Assume $\depth(R)\geq 1$ and fix an $R$-regular element $x\in\m$.
The exact sequence 
$\zeta=(0\to R\xra{x}R\to R/xR\to 0)$
is not $\hom_R(\catgpc,-)$-exact.
Hence, the natural map
$$\gpccomp^1(R/xR,R)\colon\gpcext^1(R/xR,R)\hookrightarrow\aext^1(R/xR,R)$$
is not surjective.  If $\dim(R)=1$, then
$\gpcext^1(R/xR,R)\not\cong\aext^1(R/xR,R)$.
\item \label{notiso01item4}
If $M$ admits a proper $\catpc$-resolution and
$\gpcpd_R(M)<\infty=\pcpd_R(M)$, then
$\gpcext^n(M,-)=0\ncong\pcext^n(M,-)$ for each $n>\gpcpd_R(M)$.
\end{enumerate}
\end{prop}

\begin{proof}
\eqref{notiso01item2} 
Since $C$ is in $\catpc$, we have
$\gpcext^n(C,k)\cong\pcext^n(C,k)=0$ for each $n\geq 1$ by 
Proposition~\ref{detect01}\eqref{detect01item1}.
On the other hand, the modules $\aext^n(C,k)$ are nonzero
because $C$ is a finitely generated module of infinite projective dimension
by Lemma~\ref{cincatbc'} using $B=R$.

\eqref{notiso01item3} 
Suppose that the sequence
$\zeta$ is $\hom_R(\catgpc,-)$-exact.
It follows that $\zeta$ is an augmented
proper $\catgpc$-resolution of $R/xR$, and so~\cite[(5.9)]{white:gpdrsm}
implies that it has 
an exact sequence of the following form 
as a summand
\begin{equation*} 
0\to C^n\to G \to R/xR\to 0.
\end{equation*}
where $n\geq 1$.  It follows that $C$ is a summand of $R$.
Because $R$ is local, this implies $C\cong R$
by~\cite[(2.4.b)]{white:gpdrsm}, a contradiction.

Theorem~\ref{extincl}\eqref{extinclitem1} now implies that
the natural inclusion
$\gpccomp^1(R/xR,R)$
is not surjective.
If $\dim(R)=1$, then $R/xR$ has finite length.  It follows  from
Lemma~\ref{useful2}\eqref{useful2item2} that
the module
$\gpcext^n(R/xR,R)$ has
finite length, as does $\aext^n(R/xR,R)$.  Because the inclusion $\gpccomp^1(R/xR,R)$
is not surjective, one has
$$\len_R(\gpcext^n(R/xR,R))<\len_R(\aext^n(R/xR,R))$$
and so
$\gpcext^n(R/xR,R)\not\cong\aext^n(R/xR,R)$.

\eqref{notiso01item4} 
The vanishing $\gpcext^n(M,-)=0$
when $n>\gpcpd_R(M)$ follows from Fact~\ref{useful}.
The nonvanishing
$\pcext^n(M,-)\neq 0$ is in~\cite[(3.2)]{takahashi:hasm}.
\end{proof}

Notice that the following lemma does not assume any relation between 
the semidualizing modules $B$ and $C$.
Also, the cases $B=R$ and $B=C$ imply $\pd_R(M)=\infty$
and $\pcpd_R(M)=\infty$.

\begin{lem} \label{notisoex}
Let $(R,\m)$ be a local ring and
let $B$ and $C$ be semidualizing $R$-modules.
Assume that there exist elements $y,z\in\m$ such that $\ann_R(y)=zR$ and
$\ann_R(z)=yR$, and set $M=C/yC$.  Then $M\in\catggpc=\catgpc\cap\catbc$ 
and so $M$ admits
a proper $\catpc$-resolution. Also,
one has $\pbpd_R(M)= \infty$. 
\end{lem}

\begin{proof}
Consider the chain complex
\begin{equation} \label{res}
\cdots\xra{y}C\xra{z}C\xra{y}C\xra{z}\cdots.
\end{equation}
We shall show that this complex is exact and that it is 
$\ahom(\catpc,-)$-exact and $\ahom(-,\catpc)$-exact.  
Once this is done, we will conclude from~\cite[(5.2)]{sather:sgc} 
that $M$ is in $\catggpc$.  Furthermore, we will know that the truncated complex
$$\cdots\xra{y}C\xra{z}C\xra{y} C\to 0$$
is a proper $\catpc$-resolution of $M$.

To see that the complex~\eqref{res} is exact, we first show
$\ann_R(zC)\subseteq yR$: 
If $w\in\ann_R(zC)$, we have $wz\in\ann_R(C)=0$
and so $w\in\ann_R(z)=yR$.
From this the obvious containment
$\ann_R(zC)\supseteq yR$ implies $\ann_R(zC)= yR$, and
by symmetry we have $\ann_R(yC)=zR$ and the desired exactness.  

To see that the complex~\eqref{res} is 
$\ahom(\catpc,-)$-exact and $\ahom(-,\catpc)$-exact, it suffices
to show that it is $\ahom(C,-)$-exact and $\ahom(-,C)$-exact.
This reduction uses
Hom-tensor adjointness on the one hand and~\cite[(1.11)]{white:gpdrsm}
on the other hand.  The isomorphism $\ahom(C,C)\cong R$ shows that
an application of either $\ahom(C,-)$ or $\ahom(-,C)$ yields the complex
\begin{equation} \label{res2}
\cdots\xra{y}R\xra{z}R\xra{y}R\xra{z}\cdots
\end{equation}
which is exact because of the assumptions $\ann_R(y)=zR$ and
$\ann_R(z)=yR$.

The fact that $M$ is in $\catgpc^f$ yields the
first two equalities in the next sequence 
$$0=\gpcpd_R(M)=\gcdim{C}_R(M)=\depth(R)-\depth_R(M)$$
while the third one
is from~\cite[(3.14)]{christensen:scatac}.
Now, suppose $\pbpd_R(M)<\infty$.  
Via the next sequence, the previous display works with~\cite{white:gpdrsm} to show 
that $M$ is in $\catpb^f$
$$\pbpd_R(M)=\depth(R)-\depth_R(M)=0.$$
This implies $0\neq M\cong B^m$ for some $m$,
and so $\ann_R(M)=\ann_R(B)=0$.  The membership $0\neq y\in\ann_R(M)$
contradicts this, and so $\pbpd_R(M)=\infty$.
\end{proof}

We follow-up with an example where the assumptions of Lemma~\ref{notisoex}
are satisfied.

\begin{ex} \label{ex}
Let $Q$ be a local ring with semidualizing module $A$.
Set $R=Q[\![X]\!]/(X^2)$ or $R=Q[\![Y,Z]\!]/(YZ)$. The $R$-module
$C=R\otimes_Q A$ is semidualizing by~\cite[(5.6)]{christensen:scatac}, and
the residues $y=\ol{X}=z$ or $y=\ol{Y}$ and $z=\ol{Z}$  satisfy the 
hypotheses of Lemma~\ref{notisoex}.
\end{ex}

We now contrast the relative cohomology theories arising from 
distinct semidualizing modules $B$ and $C$.
With Proposition~\ref{notiso01}\eqref{notiso01item2},
part~\eqref{notiso02item2} shows 
$\pcext^n(-,-)\ncong\pbext^n(-,-)\ncong\gpcext^n(-,-)$.
Part~\eqref{notiso02item3} shows
$\gpbext^n(-,-)\ncong\gpcext^n(-,-)$ 
and again
$\pbext^n(-,-)\ncong\gpcext^n(-,-)$.
Part~\eqref{notiso02item5} shows 
$\gpbext^n(-,-)\ncong\pcext^n(-,-)$.
Note that Lemmas~\ref{cincatbc} and~\ref{cincatbc'} 
contain analyses of  the conditions 
$\gpcpd_R(B)<\infty$ and $C\not\cong B$.

\begin{prop} \label{notiso02}
Let $(R,\m,k)$ be a local ring and
let $B$ and $C$ be semidualizing $R$-modules
such that $\gpcpd_R(B)<\infty$ and $C\not\cong B$.
\begin{enumerate}[\quad\rm(a)]
\item \label{notiso02item2}
If $n\geq 0$, then
$\pbext^n(C,k)\neq 0$.
\item \label{notiso02item3}
Assume $\depth(R)\geq 1$ and fix an $R$-regular element $x\in\m$.
The exact sequence
$\zeta=(0\to B\xra{x}B\to B/xB\to 0)$
is not $\ahom(\catgpc,-)$-exact, and so
the natural inclusions
\begin{align*}
\gpcpbcomp^1(B/xB,B)&\colon\gpcext^1(B/xB,B)\hookrightarrow\pbext^1(B/xB,B)\\
\gpcgpbcomp^1(B/xB,B)&\colon\gpcext^1(B/xB,B)\hookrightarrow\gpbext^1(B/xB,B)
\end{align*}
are not onto.  If $\dim(R)=1$, then
$\gpcext^1(B/xB,B)\not\cong\pbext^1(B/xB,B)$
and $\gpcext^1(B/xB,B)\not\cong\gpbext^1(B/xB,B)$.
\item \label{notiso02item5}
If $C$ admits a proper $\catgpb$-resolution,
then $\pcext^n(C,-)=0\neq\gpbext^n(C,-)$ for each $n\geq 1$.
\end{enumerate}
\end{prop}

\begin{proof}
\eqref{notiso02item2} 
As $\pd_R(\hom_R(B,C))=\infty$ by Lemma~\ref{cincatbc'}, 
the $n$th Betti number
$\beta^R_n(\ahom(B,C))$ is nonzero for each $n\geq 0$.
Using Fact~\ref{relabs0}\eqref{relabs0item1}, the membership $C\in\catbb$
from Lemma~\ref{cincatbc}
yields the first isomorphism in the following sequence
\begin{align*}
\pbext^n(C,k)
&\cong\aext^n(\ahom(B,C),\ahom(B,k)) \\
&\cong\aext^n(\ahom(B,C),k^{\beta^R_0(B)}) 
\cong k^{\beta^R_n(\ahom(B,C))\beta^R_0(B)}
\neq 0
\end{align*}
while the others are standard.   

\eqref{notiso02item3} 
The sequence $\zeta$
is $\catpb$-proper and $\catgpb$-proper by Fact~\ref{useful}.
Suppose that $\zeta$ is $\hom_R(\catgpc,-)$-exact. 
Because $B$ is totally $C$-reflexive by Lemma~\ref{cincatbc},
this sequence is an augmented
proper $\catgpc$-resolution of $B/xB$.  From~\cite[(5.9)]{white:gpdrsm}
we know that $\zeta$ has an exact sequence of the following form 
as a direct summand
\begin{equation*} 
0\to C^n\to G \to B/xB\to 0
\end{equation*}
where $n\geq 1$.  It follows that $C$ is a summand of $B$.
Because $R$ is local, this implies $C\cong B$
by~\cite[(2.4.b)]{white:gpdrsm}, a contradiction.
The remainder of~\eqref{notiso02item3} is  verified as in
the proof of Proposition~\ref{notiso01}\eqref{notiso01item3}.

\eqref{notiso02item5} 
Since $C$ is in $\catpc$, we have
$\pcext^n(C,-)=0$ for each $n\geq 1$.
Lemma~\ref{cincatbc'} implies
$\gpbpd_R(C)=\infty$;
arguing as in~\cite[(4.2.2.a)]{avramov:aratc}, we conclude
$\gpbext^n(C,-)\neq 0$ for each $n\geq 0$.
\end{proof}

\begin{disc} \label{resdisc}
In light of the hypothesis ``$C$ admits a proper $\catgpb$-resolution''
in Proposition~\ref{notiso02}\eqref{notiso02item5},
we note that this condition is satisfied when $R$ admits a dualizing complex
and $B=R$
by~\cite[(2.11)]{jorgensen:egprtc}.
As of the writing of this paper, the authors do not know if this condition
holds in general.
\end{disc}

We conclude this paper with 
dual versions of the above results in this section.

\begin{prop} \label{notiso01'}
Let $(R,\m,k)$ be a local ring and
$C$ a semidualizing $R$-module
such that $C\not\cong R$. Let $E$ denote the $R$-injective hull of $k$.
\begin{enumerate}[\quad\rm(a)]
\item \label{notiso01'item2}
If $n\geq 1$, then
$\icext^n(-,\ahom(C,E))=0=\gicext^n(-,\ahom(C,E))$ and
$\aext^n(-,\ahom(C,E))\neq 0$.
\item \label{notiso01'item3}
Assume that $R$ is complete and $\depth(R)\geq 1$.
Fix an $R$-regular element $x\in\m$, and
set $K=\ker(E\xra{x}E)\cong\ahom(R/xR,E)$.
The exact sequence $\zeta=(0\to K\xra{\iota} E\xra{x}E\to 0)$
is not $\ahom(-,\catgic)$-exact, and so
the map
$$\giccomp^1(E,K)\colon\gicext^1(E,K)\hookrightarrow\aext^1(E,K)$$
is not surjective. 
\item \label{notiso01'item4}
If $M$ admits a proper $\catic$-coresolution and
$\gicid_R(M)<\infty=\icid_R(M)$, then
$\gicext^n(-,M)=0\neq\icext^n(-,M)$ for each $n>\gicid_R(M)$.
\end{enumerate}
\end{prop}

\begin{proof}
The proofs of~\eqref{notiso01'item2} 
and~\eqref{notiso01'item4} are dual to the parallel parts
of Proposition~\ref{notiso01}.

\eqref{notiso01'item3}
Because $E$ is injective, it is divisible, so the sequence $\zeta$ is exact.
As in the proof of Proposition~\ref{notiso01}\eqref{notiso01item3}
it suffices to show that $\zeta$
is not $\ahom(-,\catgic)$-exact.  

Because $R$ is complete, there is an isomorphism $\ahom(E,E)\cong R$.
Applying the exact functor $\ahom(-,E)$ to $\zeta$
yields the exact sequence
\begin{equation} \label{notprope}
0\to R\xra{x}R\xra{\pi} R/xR\to 0.
\end{equation}
Proposition~\ref{notiso01}\eqref{notiso01item3} shows that
there exists a module $G\in\catgpc$ such that the following sequence is not exact
\begin{equation} \label{notexact}
0\to \ahom(G,R)\xra{x}\ahom(G,R)\xra{\ahom(G,\pi)}\ahom(G,R/xR)\to 0.
\end{equation}
Lemma~\ref{gflat}\eqref{gflatitem2}
implies that the $R$-module $G^{\vee}=\ahom(G,E)$ is in $\catgic$.
To complete the proof, we  show that the complex
$$0\to \ahom(E,G^{\vee})\xra{x}\ahom(E,G^{\vee})
\xra{\ahom(\iota,G^{\vee})} \ahom(K,G^{\vee})\to 0$$
is not exact. The ``swap'' isomorphism
$\ahom(-,G^{\vee})\cong\ahom(G,(-)^{\vee})$ shows that
this sequence is isomorphic to~\eqref{notexact}, which is not exact.
\end{proof}

\begin{lem} \label{notisoex'}
Let $(R,\m)$ be a local ring and
$C$  a semidualizing $R$-module.
Let $E$ denote the $R$-injective hull of $k$.
Assume that there exist elements $y,z\in\m$ such that $\ann_R(y)=zR$ and
$\ann_R(z)=yR$, and set $M=\ahom(C,E)/y\ahom(C,E)$.  
Then $M\in\catggic=\catgic\cap\catac$ 
and so $M$ admits
a proper $\catic$-coresolution. Also,
one has $\icid_R(M)= \infty$. 
\end{lem}

\begin{proof}
The isomorphisms $\ahom(E,E)\cong\comp{R}$ and
$C\otimes_R\ahom(C,E)\cong E$ yield the following containments
$$
0
\subseteq\ann_R(\ahom(C,E))
\subseteq\ann_R(E)
\subseteq\ann_R(\comp{R})
=0$$
and so $\ann_R(\ahom(C,E))=0$. 
As in the proof of Lemma~\ref{notisoex},
it follows that the following complex 
is exact 
\begin{equation} \label{compres}
\cdots\xra{y}\ahom(C,E)\xra{z}\ahom(C,E)\xra{y}\ahom(C,E)\xra{z}\cdots.
\end{equation}
We claim that this complex is
also $\ahom(-,\catic)$-exact and $\ahom(\catic,-)$-exact.
From this it will follow that the complex
\begin{equation} \label{iccores}
0\to\ahom(C,E)\xra{z}\ahom(C,E)\xra{y}\ahom(C,E)\xra{z}\cdots
\end{equation}
is a proper $\catic$-coresolution of $M$ and $M\in\catggic=\catgic\cap\catac$.
For  injective $R$-modules $I$ and $J$, the module
$\ahom(I,J)$ is $R$-flat. The first isomorphism in the following sequence
is Hom-tensor adjunction
$$\ahom(\ahom(C,I),\ahom(C,J))
\cong\ahom(C\otimes_R\ahom(C,I),J)
\cong\ahom(I,J)$$ 
and the second  one follows from the membership $I\in\catbc$.
Let $X$ denote the exact complex~\eqref{res2} from the proof of Lemma~\ref{notisoex}.
The displayed isomorphisms show that
an application of the functor $\ahom(\ahom(C,I),-)$ 
to the complex~\eqref{compres} yields the complex
$X\otimes_R\ahom(I,E)$.  As $X$ is exact and $\ahom(I,E)$ is flat,
the complex $X\otimes_R\ahom(I,E)$ is exact, and so~\eqref{compres} is $\ahom(\catic,-)$-exact.
Similarly, it is $\ahom(-,\catic)$-exact, as desired.

We conclude by showing $\icid_R(M)= \infty$. 
Applying the functor
$C\otimes_R -$ to the complex~\eqref{iccores} yields
an injective resolution of $C\otimes_R M$
$$0\to E \xra{z} E\xra{y} E\xra{z}\cdots.$$
This uses the memberships $M,\ahom(C,E)\in\catac$.
The fact that this resolution is minimal and nonterminating provides
the first equality in the following sequence
$$\infty=\id_R(C\otimes_R M)=\icid_R(M)$$
while the second equality is from
Fact~\ref{dimension}\eqref{dimensionitem2}.
\end{proof}

\begin{prop} \label{notiso02'}
Let $(R,\m,k)$ be a local ring and
let $B$ and $C$ be semidualizing $R$-modules
such that $\gpcpd_R(B)<\infty$ and $C\not\cong B$.
Let $E$ denote the $R$-injective hull of $k$, and set
$\mat{(-)}=\ahom(-,E)$.
\begin{enumerate}[\quad\rm(a)]
\item \label{notiso02'item2}
If $n\geq 0$, then
$\ibext^n(-,\mat{C})\neq 0$.
\item \label{notiso02'item3}
Assume that $R$ is complete and $\depth(R)\geq 1$.
Fix an $R$-regular element $x\in\m$, and
set $K= \mat{(B/xB)}$.
The sequence $\zeta=(0\to K\to \mat{B}\xra{x}\mat{B}\to 0)$
is exact but not $\ahom(-,\catgic)$-exact, and so
the natural inclusions
\begin{align*}
\gicibcomp^1(K,\mat{B})&\colon\gicext^1(K,\mat{B})\hookrightarrow\ibext^1(K,\mat{B})\\
\gicgibcomp^1(K,\mat{B})&\colon\gicext^1(K,\mat{B})\hookrightarrow\gibext^1(K,\mat{B})
\end{align*}
are not surjective. 
\item \label{notiso02'item5}
One has $\icext^n(-,\mat{C})=0\neq\gibext^n(-,\mat{C})$ for each $n\geq 1$.
\end{enumerate}
\end{prop}

\begin{proof}
\eqref{notiso02'item2}
This follows from Lemma~\ref{cincatbc'}
using~\cite[(3.2.b)]{takahashi:hasm}.

\eqref{notiso02'item3}
Consider the exact sequence
\begin{equation} \label{notprope'}
0\to B\xra{x}B\xra{\pi} B/xB\to 0.
\end{equation}
Apply the exact functor $\mat{(-)}$ to 
show that $\zeta$ is exact.
The module $\mat{B}$ is in $\catib$, so $\zeta$ is
an augmented $\catib$-coresolution of $K$.  
As in the proof of Proposition~\ref{notiso01}\eqref{notiso01item3}
it suffices to show that $\zeta$
is not $\ahom(-,\catgic)$-exact.  

Proposition~\ref{notiso02}\eqref{notiso02item3} shows that
there
exists a module $G\in\catgpc$ such that the following sequence is not exact.
\begin{equation} \label{notexact'}
0\to \ahom(G,B)\xra{x}\ahom(G,B)\xra{\ahom(G,\pi)}\ahom(G,B/xB)\to 0
\end{equation}
Lemma~\ref{gflat}\eqref{gflatitem2}
implies $G^{\vee}\in\catgic$.
To complete the proof, we  show that the complex
\begin{equation} \label{notexact''}
0\to \ahom(\mat{B},G^{\vee})\xra{x}\ahom(\mat{B},G^{\vee})
\xra{\ahom( \mat{\pi},G^{\vee})} \ahom(K,G^{\vee})\to 0
\end{equation}
is not exact. Because $R$ is complete,
the following natural  isomorphisms are valid on the
category of finitely generated $R$-modules
$$\ahom(\mat{(-)},G^{\vee})\cong\ahom(G,(-)^{\vee\vee})
\cong\ahom(G,-)$$ 
and so
the sequence~\eqref{notexact''} is isomorphic to~\eqref{notexact'}, which is not exact.

\eqref{notiso02'item5}
As in the proof of Proposition~\ref{notiso02}\eqref{notiso02item5},
it suffices to observe that 
Lemma~\ref{cincatbc'} implies
$\gibid_R(\mat{C})=\infty$.  
\end{proof}

\section*{Acknowledgments}

The authors are grateful 
to Hans-Bj\o rn Foxby, Anders Frankild and Amelia Taylor for sharing
Lemma~\ref{cincatbc}
and
to Ryo Takahashi for sharing Example~\ref{ctrex2}.

\providecommand{\bysame}{\leavevmode\hbox to3em{\hrulefill}\thinspace}
\providecommand{\MR}{\relax\ifhmode\unskip\space\fi MR }
\providecommand{\MRhref}[2]{%
  \href{http://www.ams.org/mathscinet-getitem?mr=#1}{#2}
}
\providecommand{\href}[2]{#2}

\end{document}